\documentclass{article}

\usepackage{amsmath} 
\usepackage{amssymb} 
\usepackage{amsthm} 
\usepackage{mathrsfs} 
\usepackage{mathtools} 
\usepackage{blkarray}
\usepackage{bm}

\usepackage{graphicx} 
\usepackage{float} 
\usepackage{algorithm}
\usepackage{algorithmic}

\usepackage{enumerate} 
\usepackage{booktabs} 
\usepackage{multirow} 
\usepackage{diagbox} 
\usepackage{subcaption} 
\usepackage{rotating} 
\usepackage{cite}

\usepackage{xcolor} 
\usepackage{hyperref} 
\usepackage{abstract} 
\usepackage{appendix} 
\usepackage{geometry} 

\newlength{\directplotheight}

\DeclareFontFamily{U}{mathx}{}
\DeclareFontShape{U}{mathx}{m}{n}{<-> mathx10}{}
\DeclareSymbolFont{mathx}{U}{mathx}{m}{n}
\DeclareMathAccent{\widecheck}{0}{mathx}{"71}
\makeatother

\newcommand*{\herm}{*}
\newcommand*{\pinv}{\dagger}

\newcommand*{\complex}{\mathbb{C}}

\newcommand*{\Mat}[1]{\bm{#1}}
\newcommand*{\tMat}[1]{\widetilde{\Mat{#1}}}
\newcommand*{\hMat}[1]{\widehat{\Mat{#1}}}
\newcommand*{\cMat}[1]{\widecheck{\Mat{#1}}}

\newcommand*{\bigO}{\mathscr{O}}
\newcommand*{\e}{\mathrm{e}}
\newcommand*{\imag}{\imath}
\newcommand*{\compl}{\mathrm{c}}

\newcommand*{\indexset}[1]{\mathcal{#1}}

\newcommand*{\mylog}[2][]{\log^{#1} #2}
\newcommand*{\myDim}{\mathrm{d}}
\newcommand*{\diag}{\operatorname{diag}}

\newcommand*{\floor}{\operatorname{floor}}

\newcommand*{\construct}{\mathrm{c}}
\newcommand*{\apply}{\mathrm{a}}
\newcommand*{\factor}{\mathrm{f}}
\newcommand*{\solve}{\mathrm{s}}
\newcommand*{\direct}{\mathrm{direct}}
\newcommand*{\pre}{\mathrm{offline}}
\newcommand*{\iter}{\mathrm{iter}}
\newcommand*{\piter}{\mathrm{piter}}

\newcommand*{\tree}[1]{\mathsf{#1}}
\newcommand*{\rootnode}{\rho}

\newcommand*{\ch}{\operatorname{ch}}

\newcommand*{\level}{\operatorname{level}}
\newcommand*{\midlevel}{h}
\newcommand*{\maxlevel}{L}
\newcommand*{\mylevel}{\ell}
\newcommand*{\mysublevel}{j}

\newcommand*{\maxleveloffset}{3}

\newcommand*{\hierarchical}{\mathcal{H}}

\newcommand*{\mydim}{d}
\newcommand*{\amp}{a}
\newcommand*{\phase}{\phi}
\newcommand*{\kernel}{k}
\newcommand*{\fiorightfunc}{f}
\newcommand*{\fioleftfunc}{u}
\newcommand*{\phasenonlin}{\rho}
\newcommand*{\phasecoeff}{c}
\newcommand*{\besselh}{H}

\newcommand*{\spacevar}{\Mat{x}}
\newcommand*{\freqvar}{\Mat{\xi}}
\newcommand*{\spacegrid}{X}
\newcommand*{\freqgrid}{\Xi}
\newcommand*{\spacevarsimple}{x}
\newcommand*{\freqvarsimple}{\xi}

\newcommand*{\spacemultiind}{\Mat{i}}
\newcommand*{\freqmultiind}{\Mat{j}}
\newcommand*{\spaceind}{i}
\newcommand*{\freqind}{j}

\newcommand*{\numdir}{n}
\newcommand*{\numtot}{N}

\newcommand*{\fiomat}{\Mat{K}}
\newcommand*{\fiobf}{\tMat{K}}
\newcommand*{\fioadjfio}{\Mat{G}}
\newcommand*{\fiohss}{\tMat{G}}
\newcommand*{\fiofactor}{\tMat{F}}
\newcommand*{\fioleftmat}{\Mat{u}}
\newcommand*{\fiorightmat}{\Mat{f}}
\newcommand*{\fiorightmatsol}{\hMat{f}}

\newcommand*{\spacetree}{\tree{T}_{\spacegrid}}
\newcommand*{\freqtree}{\tree{T}_{\freqgrid}}

\newcommand*{\spaceindset}{\indexset{I}}
\newcommand*{\freqindset}{\indexset{J}}
\newcommand*{\spacesubind}{t}

\newcommand*{\spacecenter}{\Mat{y}}
\newcommand*{\freqcenter}{\Mat{\eta}}
\newcommand*{\spaceinterp}{\Mat{z}}
\newcommand*{\freqinterp}{\Mat{\gamma}}
\newcommand*{\spaceinterpset}{Z}
\newcommand*{\freqinterpset}{\Gamma}

\newcommand*{\spacenode}{\tau}
\newcommand*{\freqnode}{\sigma}
\newcommand*{\spacenodepar}{\alpha}
\newcommand*{\freqnodechild}{\beta}

\newcommand*{\bfnumchild}{m}

\newcommand*{\bfrank}{r}
\newcommand*{\fiores}{\delta}

\newcommand*{\bfinterpfuncleft}{p}
\newcommand*{\bfinterpfuncright}{q}

\newcommand*{\BF}{\mathrm{BF}}
\newcommand*{\BFM}{\Mat{M}}
\newcommand*{\BFU}{\Mat{U}}
\newcommand*{\BFu}{u}
\newcommand*{\BFV}{\Mat{V}}
\newcommand*{\BFv}{v}
\newcommand*{\BFG}{\Mat{B}}
\newcommand*{\BFg}{b}
\newcommand*{\BFH}{\Mat{C}}

\newcommand*{\HSStree}{\tree{T}}
\newcommand*{\HSS}{\mathrm{HSS}}
\newcommand*{\HSSH}{\Mat{H}}
\newcommand*{\HSSD}{\Mat{D}}
\newcommand*{\HSSU}{\Mat{U}}
\newcommand*{\HSSB}{\Mat{B}}
\newcommand*{\bg}{\mathrm{big}}
\newcommand*{\HSSPOST}{\Mat{A}}

\newcommand*{\HSSindset}{\indexset{J}}
\newcommand*{\HSSnode}{\tau}
\newcommand*{\HSSnodech}{\alpha}
\newcommand*{\HSSrownode}{\tau}
\newcommand*{\HSScolnode}{\sigma}

\newcommand*{\HSSnumtot}{\numtot}

\newcommand*{\HSSrank}{r}
\newcommand*{\HSSsample}{s}
\newcommand*{\oversampling}{p}
\newcommand*{\HSSleafsize}{m}

\newcommand*{\HSSQ}{\Mat{Q}}
\newcommand*{\HSShU}{\hMat{U}}
\newcommand*{\HSScD}{\cMat{D}}
\newcommand*{\HSSL}{\Mat{L}}
\newcommand*{\HSShD}{\hMat{D}}

\newcommand*{\HSSc}{\Mat{c}}
\newcommand*{\HSShc}{\hMat{c}}
\newcommand*{\HSSy}{\Mat{y}}
\newcommand*{\HSShy}{\hMat{y}}
\newcommand*{\HSSb}{\Mat{b}}
\newcommand*{\HSShb}{\hMat{b}}

\newcommand*{\HSSOM}{\Mat{\Omega}}
\newcommand*{\HSSY}{\Mat{Y}}
\newcommand*{\HSSGA}{\Mat{\Gamma}}

\newtheorem{theorem}{Theorem}[section]
\newtheorem{definition}[theorem]{Definition}

\numberwithin{equation}{section}
\numberwithin{figure}{section}
\numberwithin{table}{section}

\providecommand{\keywords}[1]
{
  \small	
  \textbf{Keywords} #1
}

\title{Approximate Inversion of Discrete Fourier Integral Operators via Hierarchically Semiseparable Matrices}
\author{
Yingzhou Li\thanks{School of Mathematical Sciences, Fudan University; Shanghai Key Laboratory for Contemporary Applied Mathematics, Fudan University, yingzhouli@fudan.edu.cn}, 
Jingyu Liu\thanks{School of Mathematical Sciences, Fudan University, jyliu22@m.fudan.edu.cn}}
\date{\today}

\begin{document}

\maketitle

\begin{abstract}
  This paper introduces a novel method for approximating the inverse of discrete Fourier integral operators~(FIOs).
  Given an \(\numtot \times \numtot\) matrix representation \(\fiomat\) of an FIO, the proposed algorithm consists of two stages.
  In the offline stage, we first construct a butterfly factorization~(BF) \(\fiobf\) of \(\fiomat\), which enables fast forward matrix-vector multiplication.
  We then construct a hierarchically semiseparable~(HSS) approximation \(\fiohss \approx \fioadjfio\), where \(\fioadjfio = \fiomat^{\herm} \fiomat\), using fast applications of \(\fiobf\) and \(\fiobf^{\herm}\) to random matrices.
  Finally, we apply the ULV factorization to the HSS matrix \(\fiohss\) to obtain an approximation \(\fiofactor \approx \fioadjfio^{-1}\).
  Combining these approximations yields an approximate inverse \(\fiomat^{-1} \approx \fiofactor \fiobf^{\herm}\).
  The offline stage has complexity \(\bigO(\numtot \mylog[2]{\numtot})\) for 1D problems and \(\bigO(\numtot^{1.5} \mylog{\numtot})\) for 2D problems.
  In the online stage, the proposed method approximates \(\fiomat^{-1} \fioleftmat\) for a given input vector \(\fioleftmat\) with complexity \(\bigO(\numtot \mylog{\numtot})\) for both 1D and 2D problems.
  The proposed method can be used either as a direct solver or as a preconditioner for iterative methods.
  Numerical results for 1D and 2D FIOs demonstrate the effectiveness of the proposed method.
\end{abstract}

\keywords{Fourier integral operator, butterfly factorization, hierarchically semiseparable matrix}

\section{Introduction} \label{sec:intro}

This paper considers the \(\mydim\)-dimensional \textit{discrete Fourier integral operator}~(FIO) of the form
\begin{equation} \label{eq:fio}
  \fioleftfunc(\spacevar) = \sum_{\freqvar \in \freqgrid} \amp(\spacevar, \freqvar) \e^{2 \pi \imag \phase(\spacevar, \freqvar)} \fiorightfunc(\freqvar), \quad \spacevar \in \spacegrid,
\end{equation}
where
\(\spacegrid = \{\spacevar_{\spacemultiind} = (\spaceind_{1} / \numdir, \cdots, \spaceind_{\mydim} / \numdir): 0 \leq \spaceind_{1}, \ldots, \spaceind_{\mydim} < \numdir\}\)
and
\(\freqgrid = \{\freqvar_{\freqmultiind} = (\freqind_{1}, \cdots, \freqind_{\mydim}): -\numdir / 2 \leq \freqind_{1}, \ldots, \freqind_{\mydim} < \numdir / 2\}\)
for a positive even integer \(\numdir\).
Here, \(\spacevar\) and \(\freqvar\) denote the \textit{spatial} and \textit{frequency} variables, respectively.
The \textit{amplitude function} \(\amp(\spacevar, \freqvar)\) is assumed to be smooth in both \(\spacevar\) and \(\freqvar\).
The \textit{phase function} \(\phase(\spacevar, \freqvar)\) is smooth in \((\spacevar, \freqvar)\) for \(\freqvar \neq \Mat{0}\) and is homogeneous of degree \(1\) in \(\freqvar\), that is,
\(\phase(\spacevar, \lambda \freqvar) = \lambda \phase(\spacevar, \freqvar)\) for any \(\lambda > 0\).
Although the spatial and frequency grids \(\spacegrid\) and \(\freqgrid\) are assumed to be uniform in this paper, the proposed method can be extended to general grids.

In matrix form, the FIO in~\eqref{eq:fio} can be written as
\begin{equation} \label{eq:fio_mat}
  \fioleftmat = \fiomat \fiorightmat,
\end{equation}
where \(\fiomat \in \complex^{\numtot \times \numtot}\), with \(\numtot = \numdir^{\mydim}\), is the matrix representation of the FIO, and \(\fioleftmat\) and \(\fiorightmat\) are the vectors in \(\complex^{\numtot}\) corresponding to \(\fioleftfunc\) and \(\fiorightfunc\), respectively.
In some applications, only the \textit{forward} computation of the FIO is required, namely, computing \(\fioleftmat = \fiomat \fiorightmat\) for a given \(\fiorightmat\).
In this paper, we focus on the \textit{inverse} computation, namely, recovering \(\fiorightmat\) from \(\fioleftmat = \fiomat \fiorightmat\) for a given \(\fioleftmat\).
Both the forward and inverse computations are fundamental in many applications.
We refer the reader to~\cite{Candes_Demanet_Ying_2007, Cheney_Borden_2009, Symes_2009} for applications of FIOs in radar imaging, seismic imaging, and wave propagation.

\subsection{Related Work} \label{subsec:related_work}

Many algorithms have been proposed in the literature for the fast forward computation of FIOs.
In~\cite{Candes_Demanet_Ying_2007}, the authors proposed a fast method based on partitioning the frequency domain into wedges.
The discrete FIO restricted to each wedge can be decomposed into two components.
The first component is a nonuniform Fourier transform, which can be computed efficiently by the nonuniform fast Fourier transform~(NUFFT)~\cite{Barnett_Magland_Klinteberg_2019, Potts_Steidl_Tasche_2001}.
The second component has a low-rank structure and can therefore be applied efficiently.
Together, these two components reduce the cost of the forward computation from \(\bigO(\numtot^{2})\) to \(\bigO(\numtot^{5 / 4} \mylog{\numtot})\).
Another class of algorithms is based on the \textit{butterfly factorization}~(BF)~\cite{Li_Yang_2017, Li_Yang_Martin_Ho_Ying_2015}, which can be viewed as an algebraic version of the \textit{butterfly algorithm}~\cite{Candes_Demanet_Ying_2009} and has complexity \(\bigO(\numtot \mylog{\numtot})\).
Recently, in~\cite{Kielstra_Shi_Luo_Qian_Liu_2025}, a tensor butterfly factorization is proposed for multidimensional FIOs.
By combining the BF with a tensor extension of the complementary low-rank property, the cost of the forward computation is further reduced to \(\bigO(\numtot)\).

Due to the availability of efficient forward and adjoint algorithms, one can solve the linear system~\eqref{eq:fio_mat} using Krylov methods such as the generalized minimal residual method~(GMRES)~\cite{Saad_2003}, or by applying the conjugate gradient method~(CG)~\cite{Hestenes_Stiefel_1952} to the associated normal equations.
However, convergence may be slow for ill-conditioned systems, making preconditioning necessary in practice.
Moreover, fast forward representations of FIOs, such as the BF, generally cannot be inverted directly.

When dealing with matrices arising from the discretization of integral operators, the \textit{hierarchical matrix}~(\(\hierarchical\)-matrix) is a powerful tool for matrix compression and fast numerical linear algebra operations~\cite{Borm_Grasedyck_Hackbusch_2003, Hackbusch_1999, Hackbusch_Khoromskij_2000, Hackbusch_Khoromskij_Sauter_2000}.
The \(\hierarchical\)-matrix is a data-sparse representation of a dense matrix, based on a hierarchical partitioning of the matrix and low-rank approximations of certain submatrices.
Related methods include the \textit{hierarchical interpolative factorization}~(HIF)~\cite{Ho_Ying_2016a, Ho_Ying_2016b, Li_Ying_2017} and the \textit{hierarchically semiseparable}~(HSS) matrix~\cite{Martinsson_Rokhlin_2005, Xia_Chandrasekaran_Gu_Li_2010}.
These methods provide efficient ways to perform both forward applications and inverse operations in linear or quasi-linear time.
Unfortunately, FIO matrices typically do not satisfy the low-rank admissibility conditions required by these methods.

In~\cite{Feliu_Ying_2021}, the authors proposed a method for approximating the inverse of FIOs based on the identity \(\fiomat^{-1} = (\fiomat^{\herm} \fiomat)^{-1} \fiomat^{\herm}\).
The matrix \(\fiomat\) is approximated by a BF \(\fiobf\), while the matrix \(\fioadjfio = \fiomat^{\herm} \fiomat\) is approximated by an \(\hierarchical\)-matrix \(\fiohss\).
This compression is achieved by a \textit{black-box} peeling algorithm~\cite{Lin_Lu_Ying_2011}, using fast applications of \(\fiobf\) and \(\fiobf^{\herm}\).
The \(\hierarchical\)-matrix \(\fiohss\) is then inverted by an HIF-based algorithm, yielding an approximation \(\fiofactor \approx \fioadjfio^{-1}\).
Finally, the product \(\fiofactor \fiobf^{\herm}\) is used as an approximation of \(\fiomat^{-1}\).

\subsection{Contributions} \label{subsec:contributions}

The main contribution of this paper is a novel method for the approximate inversion of FIOs.
Instead of using an \(\hierarchical\)-matrix to approximate \(\fioadjfio\) and an HIF-based algorithm for the inversion, we use an HSS matrix together with its ULV factorization~\cite{Xia_Chandrasekaran_Gu_Li_2010} to achieve the same goal.
Specifically, the proposed algorithm consists of two stages.
The offline stage consists of three steps.
First, a butterfly factorization \(\fiobf\) of \(\fiomat\) is constructed using the algorithm in~\cite{Li_Yang_2017}.
Second, an HSS approximation \(\fiohss\) of \(\fioadjfio = \fiomat^{\herm} \fiomat\) is constructed by a randomized black-box algorithm.
This algorithm modifies the method in~\cite{Levitt_Martinsson_2024} by using independent test matrices at different levels, allowing the sampling to be adapted to the numerical ranks at each level.
Third, the ULV factorization of \(\fiohss\) is computed following the algorithm in~\cite{Xia_Chandrasekaran_Gu_Li_2010}, yielding an approximation \(\fiofactor \approx \fioadjfio^{-1}\).
Combining these approximations gives an approximation of the inverse of \(\fiomat\): \(\fiomat^{-1} \approx \fiofactor \fiobf^{\herm}\).
The offline stage has complexity \(\bigO(\numtot \mylog[2]{\numtot})\) for 1D problems and \(\bigO(\numtot^{1.5} \mylog{\numtot})\) for 2D problems.
In the online stage, for a given input vector \(\fioleftmat\), the approximate solution of \(\fiorightmat = \fiomat^{-1} \fioleftmat\) is computed as
\(\fiorightmat \approx \fiofactor \fiobf^{\herm} \fioleftmat\).
The online stage has complexity \(\bigO(\numtot \mylog{\numtot})\) for both 1D and 2D problems.
The proposed method can be used either as a direct solver or as a preconditioner for iterative methods.
We tested the proposed method on several examples of 1D and 2D FIOs, including constant-amplitude and variable-amplitude FIOs from the generalized Radon transform~\cite{Candes_Demanet_Ying_2007}.
The numerical results demonstrate the effectiveness and efficiency of the proposed method.

\subsection{Organization} \label{subsec:organization}

The rest of the paper is organized as follows.
Section~\ref{sec:preliminaries} reviews the butterfly factorization, the HSS matrix, and the ULV factorization, which are the main tools used in the proposed algorithm.
Section~\ref{sec:hss_blackbox_construction} presents a new algorithm for the black-box construction of HSS matrices.
Section~\ref{sec:approximate_inversion_fio} describes the proposed method for the approximate inversion of FIOs via HSS matrices.
Section~\ref{sec:numerical_results} reports numerical results for the proposed method.
Finally, Section~\ref{sec:conclusions} concludes the paper.

\section{Preliminaries} \label{sec:preliminaries}

\subsection{Notations} \label{subsec:notations}

For a set \(\indexset{J}\), we use \(|\indexset{J}|\) to denote its cardinality.
We use MATLAB notation to denote submatrices: \(\Mat{A}(\indexset{R}, \indexset{C})\) denotes the submatrix of \(\Mat{A}\) with row indices in \(\indexset{R}\) and column indices in \(\indexset{C}\).
A colon in an index denotes all indices in that dimension; for example, \(\Mat{A}(:, \indexset{C})\) denotes the submatrix of \(\Mat{A}\) with all rows and column indices in \(\indexset{C}\).
For a matrix with a superscript \(\diamond\), we write \(\Mat{A}^{\diamond, \herm}\) to denote \((\Mat{A}^{\diamond})^{\herm}\).

Tree structures are used in both the BF and the HSS matrix.
For a node \(\tau\) in a tree \(\tree{T}\), we use \(\ch(\tau)\) to denote the set of child nodes of \(\tau\).
The level of a node in \(\tree{T}\) is defined recursively as follows.
For the root node \(\rootnode\), \(\level(\rootnode) = 0\).
If \(\HSSnode\) is a nonleaf node and \(\HSSnodech \in \ch(\HSSnode)\), then \(\level(\HSSnodech) = \level(\HSSnode) + 1\).
The maximum level of the tree is denoted by \(L\).

\subsection{Butterfly Factorization} \label{subsec:bf}

In this section, we briefly review the butterfly
factorization~(BF) for FIOs~\cite{Li_Yang_2017, Li_Yang_Martin_Ho_Ying_2015}.
Let
\begin{equation*}
\kernel(\spacevar, \freqvar)
=
\amp(\spacevar, \freqvar)
\e^{2 \pi \imag \phase(\spacevar, \freqvar)}  
\end{equation*}
be the kernel associated with \(\fiomat\).
The spatial and frequency domains are recursively partitioned into
trees \(\spacetree\) and \(\freqtree\) of maximum level \(\maxlevel = \bigO(\mylog{\numtot})\).
We assume that \(\maxlevel\) is even.
For nodes \(\spacenode \in \spacetree\) and
\(\freqnode \in \freqtree\), let
\(\spaceindset_{\spacenode}\) and \(\freqindset_{\freqnode}\)
denote their corresponding index sets.
Define the corresponding local grids by
\(\spacegrid_{\spacenode}
=
\{\spacevar_{\spacemultiind}: \spacemultiind \in \spaceindset_{\spacenode}\}\)
and
\(\freqgrid_{\freqnode}
=
\{\freqvar_{\freqmultiind}: \freqmultiind \in \freqindset_{\freqnode}\}\).
The complementary low-rank property~\cite{Li_Yang_Martin_Ho_Ying_2015} states that
\(\fiomat(\spaceindset_{\spacenode},\freqindset_{\freqnode})\)
is numerically low-rank whenever \(\level(\spacenode) = \mylevel\) and \(\level(\freqnode) = \maxlevel-\mylevel\).
This property leads to the approximation
\begin{equation} \label{eq:bf}
  \fiomat
  \approx
  \fiobf
  =
  \BFU^{[\maxlevel]}
  \BFG^{[\maxlevel]}
  \cdots
  \BFG^{[\midlevel+1]}
  \BFM^{[\midlevel]}
  \BFH^{[\midlevel+1]}
  \cdots
  \BFH^{[\maxlevel]}
  \BFV^{[\maxlevel]},
\end{equation}
where \(\midlevel=\maxlevel / 2\) is the middle level of the trees.

The factors in~\eqref{eq:bf} can be constructed using Chebyshev interpolation~\cite{Li_Yang_2017}.
Let \(\bfrank\) be the interpolation rank.
For each spatial box \(\spacenode\), let \(\spacecenter_{\spacenode}\), \(\spaceinterpset_{\spacenode} = \{\spaceinterp_{\spacenode;\spaceind}\}_{\spaceind = 1}^{\bfrank}\) and \(\{\bfinterpfuncleft_{\spacenode;\spaceind}\}\) be its center, the Chebyshev points, and
the associated Lagrange functions.
Define the center \(\freqcenter_{\freqnode}\), the Chebyshev points \(\freqinterpset_{\freqnode} = \{\freqinterp_{\freqnode;\freqind}\}_{\freqind = 1}^{\bfrank}\), and the Lagrange functions
\(\{\bfinterpfuncright_{\freqnode;\freqind}\}\)
analogously for each frequency box.

For \(\midlevel\leq\mylevel\leq\maxlevel\), the phase-modulated interpolation functions are defined by
\begin{equation*}
  \BFu^{[\mylevel]}_{\spacenode, \freqnode; \spaceind}(\spacevar) = \e^{2 \pi \imag (\phase(\spacevar, \freqcenter_{\freqnode}) - \phase(\spaceinterp_{\spacenode; \spaceind}, \freqcenter_{\freqnode}))} \bfinterpfuncleft_{\spacenode; \spaceind}(\spacevar)
  \quad \text{and} \quad
  \BFv^{[\mylevel]}_{\spacenode, \freqnode; \freqind}(\freqvar) = \e^{2 \pi \imag (\phase(\spacecenter_{\spacenode}, \freqvar) - \phase(\spacecenter_{\spacenode}, \freqinterp_{\freqnode; \freqind}))} \bfinterpfuncright_{\freqnode; \freqind}(\freqvar).
\end{equation*}
For the first definition,
\(\level(\spacenode)=\mylevel\) and
\(\level(\freqnode)=\maxlevel-\mylevel\);
for the second definition, the spatial and frequency levels are
\(\maxlevel-\mylevel\) and \(\mylevel\), respectively.

At the middle level, the phase can be decomposed as
\begin{equation*}
\phase(\spacevar,\freqvar)
=
\phase(\spacevar,\freqcenter_{\freqnode})
+
\phase(\spacecenter_{\spacenode},\freqvar)
+
\fiores_{\spacenode,\freqnode}(\spacevar,\freqvar),  
\end{equation*}
where
\(\fiores_{\spacenode,\freqnode}(\spacevar,\freqvar) = \phase(\spacevar, \freqvar) - \phase(\spacevar,\freqcenter_{\freqnode}) - \phase(\spacecenter_{\spacenode},\freqvar)\).
The function \(\e^{2 \pi \imag \fiores_{\spacenode,\freqnode}(\spacevar,\freqvar)}\) can be accurately approximated by Chebyshev
interpolation~\cite{Candes_Demanet_Ying_2009}.
Together with the smoothness of the amplitude, this gives
\begin{equation} \label{eq:bf_kernel_func_low_rank}
\kernel(\spacevar, \freqvar)
\approx
\sum_{\spaceind = 1}^{\bfrank}
\sum_{\freqind = 1}^{\bfrank}
\BFu^{[\midlevel]}_{\spacenode, \freqnode; \spaceind}(\spacevar)
\kernel(\spaceinterp_{\spacenode;\spaceind}, \freqinterp_{\freqnode;\freqind})
\BFv^{[\midlevel]}_{\spacenode, \freqnode; \freqind}(\freqvar),
\quad
\spacevar\in\spacegrid_{\spacenode},
\
\freqvar\in\freqgrid_{\freqnode}.
\end{equation}
Let
\(\BFU^{[\midlevel]}_{\spacenode, \freqnode}\)
and
\(\BFV^{[\midlevel]}_{\spacenode, \freqnode}\)
be obtained by evaluating the corresponding left and right interpolation functions on
\(\spacegrid_{\spacenode}\) and
\(\freqgrid_{\freqnode}\), respectively.
Then
\begin{equation*}
  \fiomat(\spaceindset_{\spacenode}, \freqindset_{\freqnode})
  =
  \kernel(\spacegrid_{\spacenode}, \freqgrid_{\freqnode})
  \approx
  \BFU^{[\midlevel]}_{\spacenode, \freqnode}
  \BFM^{[\midlevel]}_{\spacenode, \freqnode}
  \BFV^{[\midlevel]}_{\spacenode, \freqnode},
\end{equation*}
where \(\BFM^{[\midlevel]}_{\spacenode,\freqnode} = \kernel(\spaceinterpset_{\spacenode}, \freqinterpset_{\freqnode})\).
Assembling these local approximations gives
\begin{equation} \label{eq:bf_mid}
  \fiomat
  \approx
  \BFU^{[\midlevel]}
  \BFM^{[\midlevel]}
  \BFV^{[\midlevel]}.
\end{equation}
The diagonal blocks of \(\BFU^{[\midlevel]}\) are formed by horizontally concatenating the local left interpolation matrices over frequency nodes, while those of \(\BFV^{[\midlevel]}\) are formed by vertically concatenating the local right interpolation matrices over spatial nodes.
The nonzero blocks of the weighted permutation matrix \(\BFM^{[\midlevel]}\) are
\(\kernel(\spaceinterpset_{\spacenode}, \freqinterpset_{\freqnode})\).

The same interpolation is used to construct the transfer matrices.
For \(\mylevel = \midlevel, \ldots, \maxlevel - 1\), let \(\spacenode\) be a spatial node at level \(\mylevel + 1\) with parent \(\spacenodepar\), and let \(\freqnode\) be a frequency node at level \(\maxlevel - \mylevel - 1\) with children \(\freqnodechild_{1}, \ldots, \freqnodechild_{\bfnumchild}\).
For \(\spacevar \in \spacegrid_{\spacenode}\), define
\begin{equation} \label{eq:bf_left_func}
  \BFu^{[\mylevel]}_{\spacenode, \freqnode}(\spacevar)
  =
  \sum_{\freqnodechild \in \ch(\freqnode)}
  \sum_{\spacesubind = 1}^{\bfrank}
  \BFu^{[\mylevel]}_{\spacenodepar, \freqnodechild; \spacesubind}(\spacevar)
  =
  \sum_{\freqnodechild \in \ch(\freqnode)}
  \sum_{\spacesubind = 1}^{\bfrank}
  \e^{2 \pi \imag
  \bigl(\phase(\spacevar, \freqcenter_{\freqnodechild}) - \phase(\spaceinterp_{\spacenodepar; \spacesubind},
  \freqcenter_{\freqnodechild})\bigr)}
  \bfinterpfuncleft_{\spacenodepar; \spacesubind}(\spacevar).
\end{equation}
Applying the phase-modulated interpolation associated with the new complementary pair \((\spacenode,\freqnode)\) gives
\begin{equation} \label{eq:bf_left_func_low_rank}
  \BFu^{[\mylevel]}_{\spacenode, \freqnode}(\spacevar)
  \approx
  \sum_{\spaceind = 1}^{\bfrank}
  \BFu^{[\mylevel + 1]}_{\spacenode, \freqnode;\spaceind}(\spacevar)
  \BFg^{[\mylevel + 1]}_{\spacenode, \freqnode; \spaceind},
\end{equation}
where
\(\BFg^{[\mylevel + 1]}_{\spacenode, \freqnode; \spaceind}
=
\sum_{\freqnodechild \in \ch(\freqnode)}
\sum_{\spacesubind = 1}^{\bfrank}
\BFu^{[\mylevel]}_{\spacenodepar, \freqnodechild; \spacesubind}(\spaceinterp_{\spacenode; \spaceind})\).
Let \(\BFU^{[\mylevel]}_{\spacenode, \freqnodechild}\) be the corresponding blocks in \(\BFU^{[\mylevel]}_{\spacenodepar, \freqnodechild}\), and let \(\BFU^{[\mylevel + 1]}_{\spacenode, \freqnode}\) be the interpolation matrix corresponding to \(\{\BFu^{[\mylevel + 1]}_{\spacenode, \freqnode; \spaceind}(\spacevar)\}\).
Applying the same interpolation separately to each summand in~\eqref{eq:bf_left_func} yields
\begin{equation*}
  \BFU^{[\mylevel]}_{\spacenode, \freqnode} =
  \begin{bmatrix}
    \BFU^{[\mylevel]}_{\spacenode, \freqnodechild_{1}} & \cdots & \BFU^{[\mylevel]}_{\spacenode, \freqnodechild_{\bfnumchild}}
  \end{bmatrix}
  \approx
  \BFU^{[\mylevel + 1]}_{\spacenode, \freqnode} \BFG^{[\mylevel + 1]}_{\spacenode, \freqnode},
\end{equation*}
where \(\BFG^{[\mylevel + 1]}_{\spacenode, \freqnode}\)
collects the interpolation coefficients obtained by evaluating
\(\BFu^{[\mylevel]}_{\spacenodepar, \freqnodechild; \spacesubind}\)
at the new Chebyshev points
\(\spaceinterp_{\spacenode; \spaceind}\), for
\(\freqnodechild \in \ch(\freqnode)\) and
\(1 \leq \spacesubind,\spaceind \leq \bfrank\).
The right transfer matrices are obtained similarly by exchanging
the roles of the spatial and frequency variables.
Assembling the local transfer matrices gives
\begin{equation} \label{eq:bf_recursive}
  \BFU^{[\mylevel]}
  \approx
  \BFU^{[\mylevel + 1]}
  \BFG^{[\mylevel + 1]},
  \quad
  \BFV^{[\mylevel]}
  \approx
  \BFH^{[\mylevel + 1]}
  \BFV^{[\mylevel + 1]},
  \quad
  \mylevel=\midlevel,\ldots,\maxlevel-1.
\end{equation}
Substituting~\eqref{eq:bf_recursive}
into~\eqref{eq:bf_mid} yields~\eqref{eq:bf}.

Both the construction and application of the BF have complexity
\(\bigO(\numtot\mylog{\numtot})\).
The adjoint has the same application cost and is obtained by
reversing the order of the factors in~\eqref{eq:bf} and taking
their conjugate transposes.

\subsection{HSS Matrix} \label{subsec:hss}

We review the definition of hierarchically semiseparable~(HSS) matrices in this section and their related algorithms in the next section.
Since the HSS matrix in this paper is used to approximate \(\fioadjfio = \fiomat^{\herm} \fiomat\), which is Hermitian positive definite~(HPD), we focus on HPD HSS matrices.

\begin{definition}[HSS tree,~\cite{Xia_Chandrasekaran_Gu_Li_2010}] \label{def:hss_tree}
    A tree \(\HSStree\) is called an \textit{HSS tree} with index set \(\HSSindset\) if each node \(\HSSnode\) of \(\HSStree\) is associated with an index set \(\HSSindset_{\HSSnode}\) satisfying the following conditions:
    \begin{enumerate}[(1)]
        \item \(\HSSindset_{\rootnode} = \HSSindset\) for the root node \(\rootnode\).
        \item For a nonleaf node \(\HSSnode\), we have
        \(\HSSindset_{\HSSnode} = \sqcup_{\HSSnodech \in \ch(\HSSnode)} \HSSindset_{\HSSnodech}\),
        where \(\sqcup\) denotes the disjoint union.
    \end{enumerate}
    We use \(\HSSnumtot = |\HSSindset|\) and \(\HSSnumtot_{\HSSnode} = |\HSSindset_{\HSSnode}|\) to denote the sizes of the index sets \(\HSSindset\) and \(\HSSindset_{\HSSnode}\), respectively.
\end{definition}

For simplicity, throughout this paper, we assume that the HSS tree is a full binary tree; that is, each nonleaf node has exactly two children, and all leaf nodes are at the same level \(\maxlevel\).
This assumption is natural for 1D problems.
For 2D problems, the HSS tree is typically a quadtree, and all algorithms and results can be extended to this setting in a straightforward manner.

\begin{definition}[HSS matrix,~\cite{Xia_Chandrasekaran_Gu_Li_2010}] \label{def:hss_matrix}
  Let \(\HSStree\) be an HSS tree with index set \(\HSSindset\).
  A matrix \(\HSSH \in \complex^{|\HSSindset| \times |\HSSindset|}\) is called an \textit{HSS matrix} associated with \(\HSStree\) if there exist matrices \(\HSSD_{\HSSnode}\), \(\HSSU_{\HSSnode}\), and \(\HSSB_{\HSSrownode, \HSScolnode}\), called HSS generators, associated with the nodes of \(\HSStree\) and satisfying the following conditions:
  \begin{enumerate}[(1)]
    \item
      For a leaf node \(\HSSnode\), \(\HSSD_{\HSSnode} = \HSSH({\HSSindset_{\HSSnode}, \HSSindset_{\HSSnode}}) \in \complex^{\HSSnumtot_{\HSSnode} \times \HSSnumtot_{\HSSnode}}\) is a dense Hermitian matrix, and \(\HSSU_{\HSSnode}^{\bg} = \HSSU_{\HSSnode} \in \complex^{\HSSnumtot_{\HSSnode} \times \HSSrank_{\HSSnode}}\) is a dense matrix.
    \item
      For a nonleaf node \(\HSSnode\) with children \(\HSSnodech_{1}\) and \(\HSSnodech_{2}\),
      \begin{equation*}
          \HSSD_{\HSSnode} 
          = \HSSH({\HSSindset_{\HSSnode}, \HSSindset_{\HSSnode}})
          =
          \begin{bmatrix}
              \HSSD_{\HSSnodech_{1}} & \HSSU_{\HSSnodech_{1}}^{\bg} \HSSB_{\HSSnodech_{2}, \HSSnodech_{1}}^{\herm} \HSSU_{\HSSnodech_{2}}^{\bg, \herm} \\
              \HSSU_{\HSSnodech_{2}}^{\bg} \HSSB_{\HSSnodech_{2}, \HSSnodech_{1}} \HSSU_{\HSSnodech_{1}}^{\bg, \herm} & \HSSD_{\HSSnodech_{2}}
          \end{bmatrix}.
      \end{equation*}
      For a nonleaf, nonroot node \(\HSSnode\) with children \(\HSSnodech_{1}\) and \(\HSSnodech_{2}\), the matrix \(\HSSU_{\HSSnode}^{\bg}\) satisfies the recursion
      \begin{equation*}
          \HSSU_{\HSSnode}^{\bg} =
          \begin{bmatrix}
              \HSSU_{\HSSnodech_{1}}^{\bg} & \\
              & \HSSU_{\HSSnodech_{2}}^{\bg}
          \end{bmatrix}
          \HSSU_{\HSSnode} \in \complex^{\HSSnumtot_{\HSSnode} \times \HSSrank_{\HSSnode}},
      \end{equation*}
      where \(\HSSU_{\HSSnode}\) is a dense matrix of size \((\HSSrank_{\HSSnodech_{1}} + \HSSrank_{\HSSnodech_{2}}) \times \HSSrank_{\HSSnode}\).
  \end{enumerate}
\end{definition}

By definition, \(\HSSH\) is a Hermitian matrix.
Furthermore, for any nonroot node \(\HSSnode\), the matrix \(\HSSU_{\HSSnode}^{\bg}\) forms a basis for the off-diagonal block \(\HSSH(\HSSindset_{\HSSnode}, \HSSindset_{\HSSnode}^{\compl})\), where \(\HSSindset_{\HSSnode}^{\compl} = \HSSindset \backslash \HSSindset_{\HSSnode}\).
This is often referred to as the \textit{shared basis property}, and \(\HSSH(\HSSindset_{\HSSnode}, \HSSindset_{\HSSnode}^{\compl})\) is called the \textit{HSS block} associated with \(\HSSnode\).
The \textit{HSS rank} of \(\HSSH\) is defined as the maximum rank among all HSS blocks, that is,
\(\HSSrank = \max_{\HSSnode \neq \rootnode} \HSSrank_{\HSSnode}\).
When the context is clear, we omit the index sets and use node labels directly to denote the corresponding blocks.
For example, \(\HSSH_{\HSSnode, \HSSnode^{\compl}} = \HSSH(\HSSindset_{\HSSnode}, \HSSindset_{\HSSnode}^{\compl})\).

The hierarchical structure of an HSS matrix is often represented in an alternative form known as the \textit{telescoping factorization}~\cite{Martinsson_2019}.
Specifically, define the block diagonal matrices
\begin{equation} \label{eq:hss_telescoping_factorization_factors}
  \begin{aligned}
    \HSSU^{[\mylevel]} & = \diag\bigl(\HSSU_{\HSSnode}: \level(\HSSnode) = \mylevel\bigr),  \quad 1 \leq \mylevel \leq \maxlevel, \\
    \HSSB^{[\mylevel]} & = \diag\bigl(\HSSB_{\HSSnode}: \level(\HSSnode) = \mylevel\bigr), \quad
    \HSSB_{\HSSnode} =
    \begin{bmatrix}
        \Mat{0} & \HSSB_{\HSSnodech_{2}, \HSSnodech_{1}}^{\herm} \\
        \HSSB_{\HSSnodech_{2}, \HSSnodech_{1}} & \Mat{0}
    \end{bmatrix}, \quad 0 \leq \mylevel \leq \maxlevel - 1,  \\
    \HSSD^{[\maxlevel]} & = \diag\bigl(\HSSD_{\HSSnode}: \level(\HSSnode) = \maxlevel\bigr).
  \end{aligned}
\end{equation}
Then the HSS matrix \(\HSSH\) admits the recursive factorization
\begin{equation} \label{eq:hss_telescoping_factorization}
    \begin{aligned}
        \HSSH^{[\maxlevel]} & = \HSSU^{[\maxlevel]} \HSSH^{[\maxlevel - 1]} \HSSU^{[\maxlevel], \herm} + \HSSD^{[\maxlevel]}, \\
        \HSSH^{[\mylevel]} & = \HSSU^{[\mylevel]} \HSSH^{[\mylevel - 1]} \HSSU^{[\mylevel], \herm} + \HSSB^{[\mylevel]}, \quad 1 \leq \mylevel \leq \maxlevel - 1, \\
        \HSSH^{[0]} & = \HSSB^{[0]},
    \end{aligned}
\end{equation}
where \(\HSSH^{[\maxlevel]} = \HSSH\).

\subsection{ULV Factorization of the HSS Matrix} \label{subsec:ulv}

Given the HSS representation of an HPD matrix, the HSS matrix can be inverted efficiently by the ULV factorization~\cite{Xia_Chandrasekaran_Gu_Li_2010}.
The ULV factorization can be viewed as a generalization of the Cholesky factorization.
Using the notation in Definition~\ref{def:hss_matrix}, the ULV factorization is performed in a bottom-up manner.
For a leaf node \(\HSSnode\), suppose that \(\HSSD_{\HSSnode} \in \complex^{\HSSnumtot_{\HSSnode} \times \HSSnumtot_{\HSSnode}}\) and \(\HSSU_{\HSSnode} \in \complex^{\HSSnumtot_{\HSSnode} \times \HSSrank_{\HSSnode}}\).
The first step, called \textit{basis elimination}, is to compute a unitary matrix \(\HSSQ_{\HSSnode} \in \complex^{\HSSnumtot_{\HSSnode} \times \HSSnumtot_{\HSSnode}}\) such that
\begin{equation} \label{eq:ulv_basis_elimination}
  \HSSU_{\HSSnode} = 
  \HSSQ_{\HSSnode}
  \begin{bmatrix}
    \Mat{0} \\
    \HSShU_{\HSSnode; 2}
  \end{bmatrix},
\end{equation}
where \(\HSShU_{\HSSnode; 2} \in \complex^{\HSSrank_{\HSSnode} \times \HSSrank_{\HSSnode}}\) is a dense matrix.
This can be achieved by applying the QL factorization to \(\HSSU_{\HSSnode}\).
The diagonal block \(\HSSD_{\HSSnode}\) is then transformed as \(\HSScD_{\HSSnode} = \HSSQ_{\HSSnode}^{\herm} \HSSD_{\HSSnode} \HSSQ_{\HSSnode}\).

Partition \(\HSScD_{\HSSnode}\) conformally with the partition of \(\HSSU_{\HSSnode}\) in~\eqref{eq:ulv_basis_elimination} as
\begin{equation} \label{eq:ulv_D_partition}
  \HSScD_{\HSSnode} = 
  \begin{bmatrix}
    \HSScD_{\HSSnode; 1, 1} & \HSScD_{\HSSnode; 2, 1}^{\herm} \\
    \HSScD_{\HSSnode; 2, 1} & \HSScD_{\HSSnode; 2, 2}
  \end{bmatrix}.
\end{equation}
A \textit{partial factorization} is then performed to eliminate the \((1, 1)\) block of \(\HSScD_{\HSSnode}\).
More precisely, the Cholesky factorization of \(\HSScD_{\HSSnode; 1, 1}\) is computed as
\(\HSScD_{\HSSnode; 1, 1} = \HSSL_{\HSSnode; 1, 1} \HSSL_{\HSSnode; 1, 1}^{\herm}\).
Then~\eqref{eq:ulv_D_partition} can be decomposed as
\begin{equation} \label{eq:ulv_partial_factorization}
  \HSScD_{\HSSnode} =
  \begin{bmatrix}
    \HSSL_{\HSSnode; 1, 1}  &  \\
    \HSShD_{\HSSnode; 2, 1} & \Mat{I}
  \end{bmatrix}
  \begin{bmatrix}
    \Mat{I}  &  \\
     & \HSShD_{\HSSnode; 2, 2}
  \end{bmatrix}
  \begin{bmatrix}
    \HSSL_{\HSSnode; 1, 1}  &  \\
    \HSShD_{\HSSnode; 2, 1} & \Mat{I}
  \end{bmatrix}^{\herm},
\end{equation}
where
\(\HSShD_{\HSSnode; 2, 1} = \HSScD_{\HSSnode; 2, 1} \HSSL_{\HSSnode; 1, 1}^{-\herm}\)
and
\(\HSShD_{\HSSnode; 2, 2} = \HSScD_{\HSSnode; 2, 2} - \HSShD_{\HSSnode; 2, 1} \HSShD_{\HSSnode; 2, 1}^{\herm}\)
is the Schur complement.

After basis elimination and partial factorization have been performed for all leaf nodes, the algorithm proceeds to level \(\maxlevel - 1\).
For a nonleaf node \(\HSSnode\) at level \(\maxlevel - 1\) with children \(\HSSnodech_{1}\) and \(\HSSnodech_{2}\), a \textit{merge} step is performed to update the HSS generators of \(\HSSnode\) as
\begin{equation} \label{eq:ulv_merge}
  \HSSD_{\HSSnode}^{+} = 
      \begin{bmatrix}
          \HSShD_{\HSSnodech_{1}; 2, 2} & \HSShU_{\HSSnodech_{1}; 2} \HSSB_{\HSSnodech_{2}, \HSSnodech_{1}}^{\herm} \HSShU_{\HSSnodech_{2}; 2}^{\herm} \\
          \HSShU_{\HSSnodech_{2}; 2} \HSSB_{\HSSnodech_{2}, \HSSnodech_{1}} \HSShU_{\HSSnodech_{1}; 2}^{\herm} & \HSShD_{\HSSnodech_{2}; 2, 2}
      \end{bmatrix}, \quad 
      \HSSU_{\HSSnode}^{+} = 
      \begin{bmatrix}
          \HSShU_{\HSSnodech_{1}; 2} & \\
          & \HSShU_{\HSSnodech_{2}; 2}
      \end{bmatrix}
      \HSSU_{\HSSnode}.
\end{equation}
The children \(\HSSnodech_{1}\) and \(\HSSnodech_{2}\) can then be ``removed'', and the node \(\HSSnode\) can be treated as a leaf node with the updated HSS generators \(\HSSD_{\HSSnode}^{+}\) and \(\HSSU_{\HSSnode}^{+}\).
These steps are repeated from bottom to top, with \(\HSSD_{\HSSnode}\) and \(\HSSU_{\HSSnode}\) replaced by \(\HSSD_{\HSSnode}^{+}\) and \(\HSSU_{\HSSnode}^{+}\), until the root node \(\rootnode\) is reached.
At the root node, the Cholesky factorization \(\HSSD_{\rootnode}^{+}= \HSSL_{\rootnode} \HSSL_{\rootnode}^{\herm}\) is computed directly.

Once the ULV factorization has been computed, the linear system \(\HSSH \HSSc = \HSSy\) can be solved efficiently by a bottom-up pass followed by a top-down pass.
First, the vector \(\HSSy\) is partitioned and assigned to the leaf nodes of the HSS tree.
Starting from each leaf node \(\HSSnode\), the vector \(\HSSy_{\HSSnode}\) is updated by
\begin{equation*}
  \HSShy_{\HSSnode} = \HSSQ_{\HSSnode}^{\herm} \HSSy_{\HSSnode} =
  \begin{bmatrix}
    \HSShy_{\HSSnode; 1} \\
    \HSShy_{\HSSnode; 2}
  \end{bmatrix}.
\end{equation*}
Then the vector \(\HSShb_{\HSSnode; 1}\) is computed by solving \(\HSSL_{\HSSnode; 1, 1} \HSShb_{\HSSnode; 1} = \HSShy_{\HSSnode; 1}\), and the vector \(\HSShy_{\HSSnode; 2}\) is updated as \(\HSShy_{\HSSnode; 2} \gets \HSShy_{\HSSnode; 2} - \HSShD_{\HSSnode; 2, 1} \HSShb_{\HSSnode; 1}\).
Next, for a nonleaf node \(\HSSnode\) at level \(\maxlevel - 1\) with children \(\HSSnodech_{1}\) and \(\HSSnodech_{2}\), the vector \(\HSSy_{\HSSnode}^{+}\) is obtained by merging \(\HSShy_{\HSSnodech_{1}; 2}\) and \(\HSShy_{\HSSnodech_{2}; 2}\) as
\begin{equation*}
  \HSSy_{\HSSnode}^{+} =
\begin{bmatrix}
  \HSShy_{\HSSnodech_{1}; 2} \\
  \HSShy_{\HSSnodech_{2}; 2}
\end{bmatrix}.
\end{equation*}
This process is repeated from bottom to top, with \(\HSSy_{\HSSnode}\) replaced by \(\HSSy_{\HSSnode}^{+}\), until the root node \(\rootnode\) is reached.
At the root node, the vectors \(\HSSb_{\rootnode}\) and
\(\HSSc_{\rootnode}\) are obtained by solving the triangular systems
\(\HSSL_{\rootnode}\HSSb_{\rootnode}
= \HSSy_{\rootnode}^{+}\)
and
\(\HSSL_{\rootnode}^{\herm}\HSSc_{\rootnode}
= \HSSb_{\rootnode}\),
respectively.

A top-down pass is then performed to compute the solution \(\HSSc\) from the root node to the leaf nodes.
Starting from the root node \(\rootnode\) with children \(\HSSnodech_{1}\) and \(\HSSnodech_{2}\), the vectors \(\HSShc_{\HSSnodech_{1}; 2}\) and \(\HSShc_{\HSSnodech_{2}; 2}\) are obtained from the partition
\begin{equation*}
  \HSSc_{\rootnode} =
  \begin{bmatrix}
    \HSShc_{\HSSnodech_{1}; 2} \\
    \HSShc_{\HSSnodech_{2}; 2}
  \end{bmatrix}.
\end{equation*}
Then, for \(1 \leq \mylevel \leq \maxlevel\) and every node \(\HSSnode\) at level \(\mylevel\), the vector \(\HSShc_{\HSSnode; 1}\) is obtained by solving
\begin{equation*}
  \HSSL_{\HSSnode; 1, 1}^{\herm} \HSShc_{\HSSnode; 1}
  =
  \HSShb_{\HSSnode; 1} - \HSShD_{\HSSnode; 2, 1}^{\herm} \HSShc_{\HSSnode; 2},
\end{equation*}
and the vector \(\HSSc_{\HSSnode}\) is given by
\begin{equation*}
  \HSSc_{\HSSnode}
  =
  \HSSQ_{\HSSnode}
  \begin{bmatrix}
    \HSShc_{\HSSnode; 1} \\
    \HSShc_{\HSSnode; 2}
  \end{bmatrix}.
\end{equation*}
If \(\HSSnode\) is a nonleaf node with children \(\HSSnodech_{1}\) and \(\HSSnodech_{2}\), then \(\HSShc_{\HSSnodech_{1}; 2}\) and \(\HSShc_{\HSSnodech_{2}; 2}\) are obtained from the partition of \(\HSSc_{\HSSnode}\).
Finally, the solution \(\HSSc\) is obtained by collecting the vectors \(\HSSc_{\HSSnode}\) over all leaf nodes \(\HSSnode\).

\section{A Black-box Construction of the HSS Matrix} \label{sec:hss_blackbox_construction}

In this section, we present a new algorithm for the black-box construction of a Hermitian HSS matrix.
The algorithm is based on the method proposed in~\cite{Levitt_Martinsson_2024}.
Using the notation in Definition~\ref{def:hss_matrix}, the HSS matrix \(\HSSH\) is constructed by sequentially applying \(\HSSH\) to random matrices.
Unlike the method in~\cite{Levitt_Martinsson_2024}, the proposed algorithm uses independent random matrices at each level, allowing the sampling to be adapted to the numerical ranks at that level.
The construction consists of two steps.
First, a bottom-up process is performed to compute all basis matrices \(\HSSU_{\HSSnode}\) and residual blocks \(\HSScD_{\HSSnode}\).
Next, a top-down process is used to complete the construction by assigning the matrices \(\HSSB_{\HSSnodech_{2}, \HSSnodech_{1}}\) and \(\HSSD_{\HSSnode}\).
In this section, we assume that the rank of the HSS blocks at level \(\mylevel\) is \(\HSSrank_{\mylevel}\), with \(\HSSrank_{0} = 0\), and that all leaf nodes have size \(\HSSnumtot_{\maxlevel}\).
We define \(\HSSleafsize_{\maxlevel} = \HSSnumtot_{\maxlevel}\) and \(\HSSleafsize_{\mylevel} = 2 \HSSrank_{\mylevel + 1}\) for \(0 \leq \mylevel \leq \maxlevel - 1\), and refer to \(\HSSleafsize_{\mylevel}\) as the ``effective leaf size'' at level \(\mylevel\).

\subsection{Construction on the Leaf Nodes} \label{susec:leaf_construction}

Starting from the leaf level \(\maxlevel\), let
\(\HSSsample_{\maxlevel} = \HSSrank_{\maxlevel} + \HSSleafsize_{\maxlevel} + \oversampling\),
where \(\oversampling\) is an oversampling parameter, e.g., \(\oversampling = 5\).
Let \(\HSSOM^{[\maxlevel]} \in \complex^{|\HSSindset| \times \HSSsample_{\maxlevel}}\) be a random matrix with i.i.d. Gaussian entries, and let
\(\HSSY^{[\maxlevel]} = \HSSH \HSSOM^{[\maxlevel]}\).
We call \(\HSSOM^{[\maxlevel]}\) and \(\HSSY^{[\maxlevel]}\) the \textit{test matrix} and \textit{sample matrix} at level \(\maxlevel\), respectively.
For each leaf node \(\HSSnode\) at level \(\maxlevel\), by considering the row block corresponding to \(\HSSnode\) in the equation \(\HSSY^{[\maxlevel]} = \HSSH \HSSOM^{[\maxlevel]}\), we obtain
\begin{equation} \label{eq:hss_construct_leaf_block}
    \HSSY^{[\maxlevel]}_{\HSSnode} = \HSSD_{\HSSnode} \HSSOM^{[\maxlevel]}_{\HSSnode} + \HSSH_{\HSSnode, \HSSnode^{\compl}} \HSSOM^{[\maxlevel]}_{\HSSnode^{\compl}}.
\end{equation}

Let
\(\HSSGA^{[\maxlevel]}_{\HSSnode}
\in \complex^{\HSSsample_{\maxlevel}
\times(\HSSrank_{\maxlevel}+\oversampling)}\)
be an orthonormal basis for the null space of
\(\HSSOM^{[\maxlevel]}_{\HSSnode}\), so that
\(\HSSOM^{[\maxlevel]}_{\HSSnode}
\HSSGA^{[\maxlevel]}_{\HSSnode}=\Mat{0}\).
Such a matrix can be computed by performing a complete QR
factorization with column pivoting~(QRCP) on
\(\HSSOM^{[\maxlevel],\herm}_{\HSSnode}\)
and taking the last
\(\HSSrank_{\maxlevel}+\oversampling\)
columns of the Q-factor.
Multiplying both sides of~\eqref{eq:hss_construct_leaf_block} by \(\HSSGA^{[\maxlevel]}_{\HSSnode}\) gives
\begin{equation*}
    \HSSY^{[\maxlevel]}_{\HSSnode} \HSSGA^{[\maxlevel]}_{\HSSnode} = \HSSH_{\HSSnode, \HSSnode^{\compl}} \HSSOM^{[\maxlevel]}_{\HSSnode^{\compl}} \HSSGA^{[\maxlevel]}_{\HSSnode}.
\end{equation*}
Since \(\HSSGA^{[\maxlevel]}_{\HSSnode}\) is orthonormal, the matrix \(\HSSY^{[\maxlevel]}_{\HSSnode} \HSSGA^{[\maxlevel]}_{\HSSnode}\) can be viewed as a sample matrix for \(\HSSH_{\HSSnode, \HSSnode^{\compl}}\) with test matrix \(\HSSOM^{[\maxlevel]}_{\HSSnode^{\compl}} \HSSGA^{[\maxlevel]}_{\HSSnode}\).
Consequently, \(\HSSU_{\HSSnode}\) can be constructed by computing an orthonormal basis for the column space of \(\HSSY^{[\maxlevel]}_{\HSSnode} \HSSGA^{[\maxlevel]}_{\HSSnode}\), for example, by QRCP.
In practice, the basis is truncated according to the prescribed accuracy.

Subsequently, the diagonal block \(\HSSD_{\HSSnode}\) is decomposed into two parts:
\begin{equation} \label{eq:hss_construct_D_partition}
  \begin{aligned}
    \HSSD_{\HSSnode}
    & = \HSSU_{\HSSnode} (\HSSU_{\HSSnode}^{\herm} \HSSD_{\HSSnode} \HSSU_{\HSSnode}) \HSSU_{\HSSnode}^{\herm}
    + \bigl(\HSSD_{\HSSnode} - \HSSU_{\HSSnode} \HSSU_{\HSSnode}^{\herm} \HSSD_{\HSSnode} \HSSU_{\HSSnode} \HSSU_{\HSSnode}^{\herm}\bigr) \\
    & = \HSSU_{\HSSnode} \HSShD_{\HSSnode} \HSSU_{\HSSnode}^{\herm} + \HSScD_{\HSSnode},
  \end{aligned}
\end{equation}
where \(\HSShD_{\HSSnode} = \HSSU_{\HSSnode}^{\herm} \HSSD_{\HSSnode} \HSSU_{\HSSnode}\) and
\(\HSScD_{\HSSnode} = \HSSD_{\HSSnode} - \HSSU_{\HSSnode} \HSShD_{\HSSnode} \HSSU_{\HSSnode}^{\herm}\).
The first term represents the projection of \(\HSSD_{\HSSnode}\) onto the column space of \(\HSSU_{\HSSnode}\), while the second term is the residual.
We will discuss the computation of \(\HSShD_{\HSSnode}\) later and first focus on the computation of \(\HSScD_{\HSSnode}\).
By rewriting \(\HSScD_{\HSSnode}\) as
\begin{equation*}
  \begin{aligned}
  \HSScD_{\HSSnode}
  & = \HSSD_{\HSSnode} - \HSSU_{\HSSnode} \HSSU_{\HSSnode}^{\herm} \HSSD_{\HSSnode} \HSSU_{\HSSnode} \HSSU_{\HSSnode}^{\herm} \\
  & = (\Mat{I} - \HSSU_{\HSSnode} \HSSU_{\HSSnode}^{\herm}) \HSSD_{\HSSnode}
  + \HSSU_{\HSSnode} \HSSU_{\HSSnode}^{\herm} \HSSD_{\HSSnode} (\Mat{I} - \HSSU_{\HSSnode} \HSSU_{\HSSnode}^{\herm}),
  \end{aligned}
\end{equation*}
it suffices to compute \((\Mat{I} - \HSSU_{\HSSnode} \HSSU_{\HSSnode}^{\herm}) \HSSD_{\HSSnode}\).
Multiplying both sides of~\eqref{eq:hss_construct_leaf_block} by \((\Mat{I} - \HSSU_{\HSSnode} \HSSU_{\HSSnode}^{\herm})\) and using the fact that \(\HSSU_{\HSSnode}\) is an orthonormal basis for \(\HSSH_{\HSSnode, \HSSnode^{\compl}}\), we obtain
\begin{equation*}
  (\Mat{I} - \HSSU_{\HSSnode} \HSSU_{\HSSnode}^{\herm}) \HSSD_{\HSSnode} \HSSOM^{[\maxlevel]}_{\HSSnode}
  =
  (\Mat{I} - \HSSU_{\HSSnode} \HSSU_{\HSSnode}^{\herm}) \HSSY^{[\maxlevel]}_{\HSSnode}.
\end{equation*}
Therefore,
\begin{equation*}
  (\Mat{I} - \HSSU_{\HSSnode} \HSSU_{\HSSnode}^{\herm}) \HSSD_{\HSSnode}
  =
  (\Mat{I} - \HSSU_{\HSSnode} \HSSU_{\HSSnode}^{\herm}) \HSSY^{[\maxlevel]}_{\HSSnode} \HSSOM^{[\maxlevel], \pinv}_{\HSSnode}.
\end{equation*}
Consequently, the residual block is computed as
\begin{equation*}
  \HSScD_{\HSSnode}
  =
  (\Mat{I} - \HSSU_{\HSSnode} \HSSU_{\HSSnode}^{\herm}) \HSSY^{[\maxlevel]}_{\HSSnode} \HSSOM^{[\maxlevel], \pinv}_{\HSSnode}
  +
  \HSSU_{\HSSnode} \HSSU_{\HSSnode}^{\herm}
  \bigl(
  (\Mat{I} - \HSSU_{\HSSnode} \HSSU_{\HSSnode}^{\herm}) \HSSY^{[\maxlevel]}_{\HSSnode} \HSSOM^{[\maxlevel], \pinv}_{\HSSnode}
  \bigr)^{\herm}.
\end{equation*}

\subsection{Construction on the Nonleaf Nodes} \label{subsec:nonleaf_construction}

Suppose \(1 \leq \mylevel \leq \maxlevel - 1\) and that the basis matrices for all nodes at levels \(\mylevel + 1, \dotsc, \maxlevel\) have been computed.
Define \(\HSSU^{[\mylevel + 1]} = \diag\bigl(\HSSU_{\HSSnode}: \level(\HSSnode) = \mylevel + 1\bigr)\).
By extracting the basis matrices, we obtain the following recursive decomposition of \(\HSSH^{[\mylevel + 1]}\), starting from \(\HSSH^{[\maxlevel]} = \HSSH\):
\begin{equation} \label{eq:hss_construct_telescoping}
    \HSSH^{[\mylevel + 1]} = \HSSU^{[\mylevel + 1]} \HSSH^{[\mylevel]} \HSSU^{[\mylevel + 1], \herm} + \HSScD^{[\mylevel + 1]},
\end{equation}
where \(\HSSH^{[\mylevel]} = \HSSU^{[\mylevel + 1], \herm} \HSSH^{[\mylevel + 1]} \HSSU^{[\mylevel + 1]}\) and \(\HSScD^{[\mylevel + 1]} = \diag\bigl(\HSScD_{\HSSnode}: \level(\HSSnode) = \mylevel + 1\bigr)\).
If all nodes at level \(\mylevel + 1\) are ``eliminated'', then \(\HSSH^{[\mylevel]}\) remains an HSS matrix with maximum level \(\mylevel\), but with updated generators.
Specifically, for every node \(\HSSnode\) at level \(\mylevel\) with children \(\HSSnodech_{1}\) and \(\HSSnodech_{2}\), the node \(\HSSnode\) becomes a leaf node of the reduced HSS tree, with its diagonal generator updated as
\begin{equation*}
    \HSSD_{\HSSnode}
    \gets
    \begin{bmatrix}
        \HSShD_{\HSSnodech_{1}} & \HSSB_{\HSSnodech_{2}, \HSSnodech_{1}}^{\herm} \\
        \HSSB_{\HSSnodech_{2}, \HSSnodech_{1}} & \HSShD_{\HSSnodech_{2}}
    \end{bmatrix}.
\end{equation*}

Let \(\HSSsample_{\mylevel} = \HSSrank_{\mylevel} + \HSSleafsize_{\mylevel} + \oversampling\), and let \(\HSSOM^{[\mylevel]} \in \complex^{2^{\mylevel} \HSSleafsize_{\mylevel} \times \HSSsample_{\mylevel}}\) be the test matrix at level \(\mylevel\).
The corresponding sample matrix is defined as \(\HSSY^{[\mylevel]} = \HSSH^{[\mylevel]} \HSSOM^{[\mylevel]}\).
Then
\begin{equation*}
  \HSSY^{[\mylevel]}
  = \HSSH^{[\mylevel]} \HSSOM^{[\mylevel]} \\
  = \HSSU^{[\mylevel + 1], \herm} \HSSH^{[\mylevel + 1]} \HSSU^{[\mylevel + 1]} \HSSOM^{[\mylevel]} \\
  = \Bigl(\HSSU^{[\mylevel + 1], \herm} \dotsb \HSSU^{[\maxlevel], \herm}\Bigr)
      \HSSH
      \Bigl(\HSSU^{[\maxlevel]} \dotsb \HSSU^{[\mylevel + 1]}\Bigr)
      \HSSOM^{[\mylevel]}.
\end{equation*}
That is, the sample matrix \(\HSSY^{[\mylevel]}\) can be obtained by first applying the basis matrices from level \(\mylevel + 1\) up to level \(\maxlevel\) to the test matrix \(\HSSOM^{[\mylevel]}\), then applying \(\HSSH\) to the resulting matrix, and finally applying the adjoint basis matrices from level \(\maxlevel\) down to level \(\mylevel + 1\).
Once \(\HSSY^{[\mylevel]}\) is obtained, the basis matrices and residual matrices for all nodes at level \(\mylevel\) can be computed by a procedure similar to that in Section~\ref{susec:leaf_construction}.
At the root node, we only need to solve for the updated diagonal block \(\HSSH^{[0]}\) from
\(\HSSY^{[0]} = \HSSH^{[0]} \HSSOM^{[0]}\).

\subsection{Postprocessing} \label{subsec:hss_construction_postprocessing}

After computing all matrices \(\HSSU_{\HSSnode}\) and \(\HSScD_{\HSSnode}\), a top-down procedure is performed to construct the matrices \(\HSSB\) and \(\HSSD\), as summarized in~\cite{Li_Liu_2025}.
For level \(0 \leq \mylevel \leq \maxlevel - 1\), suppose that \(\HSSnode\) is a nonleaf node at level \(\mylevel\) with children \(\HSSnodech_{1}\) and \(\HSSnodech_{2}\).
We define \(\HSSD_{\HSSnode}^{[\mylevel]}\) and \(\HSSU_{\HSSnode}^{[\mylevel + 1]}\) as the submatrices of \(\HSSD^{[\mylevel]}\) and \(\HSSU^{[\mylevel + 1]}\) corresponding to the node \(\HSSnode\), respectively.
The equation~\eqref{eq:hss_construct_telescoping} corresponding to the node \(\HSSnode\) is
\begin{equation} \label{eq:hss_construct_postprocess_decompose}
    \begin{aligned}
      \HSSH_{\HSSnode}^{[\mylevel + 1]}
      & = 
        \HSSU_{\HSSnode}^{[\mylevel + 1]} \HSSH_{\HSSnode}^{[\mylevel]} \HSSU_{\HSSnode}^{[\mylevel + 1], \herm}
        +
        \HSScD_{\HSSnode}^{[\mylevel + 1]} \\ 
      & =
        \HSSU_{\HSSnode}^{[\mylevel + 1]}
        \Bigl(
          \HSSU_{\HSSnode}^{[\mylevel]} \HSSH_{\HSSnode}^{[\mylevel - 1], +} \HSSU_{\HSSnode}^{[\mylevel], \herm}
          +
          \HSSD_{\HSSnode}^{[\mylevel]}
        \Bigr)
        \HSSU_{\HSSnode}^{[\mylevel + 1], \herm}
        +
        \HSScD_{\HSSnode}^{[\mylevel + 1]} \\ 
      & =
        \HSSU_{\HSSnode}^{[\mylevel + 1]} \HSSU_{\HSSnode}^{[\mylevel]} \HSSH_{\HSSnode}^{[\mylevel - 1], +}
        \HSSU_{\HSSnode}^{[\mylevel], \herm} \HSSU_{\HSSnode}^{[\mylevel + 1], \herm}
        +
        \HSSU_{\HSSnode}^{[\mylevel + 1]} \HSSD_{\HSSnode}^{[\mylevel]} \HSSU_{\HSSnode}^{[\mylevel + 1], \herm}
        +
        \HSScD_{\HSSnode}^{[\mylevel + 1]}.
    \end{aligned}
\end{equation}
Here, \(\HSSH_{\HSSnode}^{[\mylevel - 1], +}\) is consistent with the telescoping factorization in~\eqref{eq:hss_telescoping_factorization}.
For the root node \(\rootnode\), the matrices are initialized as
\(\HSSD_{\rootnode}^{[0]} = \HSSH^{[0]}\),
\(\HSSU_{\rootnode}^{[0]} = \Mat{I}\),
and \(\HSSH_{\rootnode}^{[- 1], +} = \Mat{0}\).
Compared with the telescoping factorization~\eqref{eq:hss_telescoping_factorization}, the only difference is that the diagonal blocks of \(\HSSD_{\HSSnode}^{[\mylevel]}\) are nonzero.
Partition \(\HSSD_{\HSSnode}^{[\mylevel]}\) as
\begin{equation*}
    \HSSD_{\HSSnode}^{[\mylevel]} = 
    \begin{bmatrix}
        \HSSPOST_{\HSSnode; \HSSnodech_{1}, \HSSnodech_{1}} & \HSSPOST_{\HSSnode; \HSSnodech_{2}, \HSSnodech_{1}}^{\herm} \\
        \HSSPOST_{\HSSnode; \HSSnodech_{2}, \HSSnodech_{1}} & \HSSPOST_{\HSSnode; \HSSnodech_{2}, \HSSnodech_{2}}
    \end{bmatrix}.
\end{equation*}
Then the second and third terms in~\eqref{eq:hss_construct_postprocess_decompose} can be rewritten as
\begin{equation*}
    \begin{aligned}
      &\HSSU_{\HSSnode}^{[\mylevel + 1]} \HSSD_{\HSSnode}^{[\mylevel]} \HSSU_{\HSSnode}^{[\mylevel + 1], \herm}
      +
      \HSScD_{\HSSnode}^{[\mylevel + 1]} \\
      & \quad =
      \begin{bmatrix}
          \HSSU_{\HSSnodech_{1}} & \\
          & \HSSU_{\HSSnodech_{2}}
      \end{bmatrix}
      \begin{bmatrix}
          \HSSPOST_{\HSSnode; \HSSnodech_{1}, \HSSnodech_{1}} & \HSSPOST_{\HSSnode; \HSSnodech_{2}, \HSSnodech_{1}}^{\herm} \\
          \HSSPOST_{\HSSnode; \HSSnodech_{2}, \HSSnodech_{1}} & \HSSPOST_{\HSSnode; \HSSnodech_{2}, \HSSnodech_{2}}
      \end{bmatrix}
      \begin{bmatrix}
          \HSSU_{\HSSnodech_{1}} & \\
          & \HSSU_{\HSSnodech_{2}}
      \end{bmatrix}^{\herm}
      +
      \begin{bmatrix}
          \HSScD_{\HSSnodech_{1}} & \\
          & \HSScD_{\HSSnodech_{2}}
      \end{bmatrix} \\
      & \quad =
      \begin{bmatrix}
          \HSSU_{\HSSnodech_{1}} & \\
          & \HSSU_{\HSSnodech_{2}}
      \end{bmatrix}
      \begin{bmatrix}
          & \HSSB_{\HSSnodech_{2}, \HSSnodech_{1}}^{\herm} \\
          \HSSB_{\HSSnodech_{2}, \HSSnodech_{1}} &
      \end{bmatrix}
      \begin{bmatrix}
          \HSSU_{\HSSnodech_{1}} & \\
          & \HSSU_{\HSSnodech_{2}}
      \end{bmatrix}^{\herm}
      +
      \begin{bmatrix}
          \HSSD_{\HSSnodech_{1}} & \\
          & \HSSD_{\HSSnodech_{2}}
      \end{bmatrix} \\
      & \quad =
      \HSSU_{\HSSnode}^{[\mylevel + 1]} \HSSB_{\HSSnode}^{[\mylevel]} \HSSU_{\HSSnode}^{[\mylevel + 1], \herm}
      +
      \HSSD_{\HSSnode}^{[\mylevel + 1]},
  \end{aligned}
\end{equation*}
where \(\HSSB_{\HSSnodech_{2}, \HSSnodech_{1}} = \HSSPOST_{\HSSnode; \HSSnodech_{2}, \HSSnodech_{1}}\) and \(\HSSD_{\HSSnodech}
=
\HSScD_{\HSSnodech}
+
\HSSU_{\HSSnodech}
\HSSPOST_{\HSSnode; \HSSnodech, \HSSnodech}
\HSSU_{\HSSnodech}^{\herm}\) for each \(\HSSnodech \in \ch(\HSSnode)\).
Therefore, equation~\eqref{eq:hss_construct_postprocess_decompose} can be reformulated as
\begin{equation*}
    \begin{aligned}
      \HSSH_{\HSSnode}^{[\mylevel + 1]} 
        & =
          \HSSU_{\HSSnode}^{[\mylevel + 1]} \HSSU_{\HSSnode}^{[\mylevel]} \HSSH_{\HSSnode}^{[\mylevel - 1], +}
          \HSSU_{\HSSnode}^{[\mylevel], \herm} \HSSU_{\HSSnode}^{[\mylevel + 1], \herm}
          +
          \HSSU_{\HSSnode}^{[\mylevel + 1]} \HSSB_{\HSSnode}^{[\mylevel]} \HSSU_{\HSSnode}^{[\mylevel + 1], \herm}
          +
          \HSSD_{\HSSnode}^{[\mylevel + 1]} \\
        & =
          \HSSU_{\HSSnode}^{[\mylevel + 1]}
          \Bigl(
            \HSSU_{\HSSnode}^{[\mylevel]} \HSSH_{\HSSnode}^{[\mylevel - 1], +} \HSSU_{\HSSnode}^{[\mylevel], \herm}
            +
            \HSSB_{\HSSnode}^{[\mylevel]}
          \Bigr)
          \HSSU_{\HSSnode}^{[\mylevel + 1], \herm}
          +
          \HSSD_{\HSSnode}^{[\mylevel + 1]} \\
        & =
          \HSSU_{\HSSnode}^{[\mylevel + 1]} \HSSH_{\HSSnode}^{[\mylevel], +} \HSSU_{\HSSnode}^{[\mylevel + 1], \herm}
          +
          \HSSD_{\HSSnode}^{[\mylevel + 1]},
    \end{aligned}
\end{equation*}
where
\(\HSSH_{\HSSnode}^{[\mylevel], +}
=
\HSSU_{\HSSnode}^{[\mylevel]} \HSSH_{\HSSnode}^{[\mylevel - 1], +} \HSSU_{\HSSnode}^{[\mylevel], \herm}
+
\HSSB_{\HSSnode}^{[\mylevel]}\)
is consistent with the telescoping factorization in~\eqref{eq:hss_telescoping_factorization}.
The process is repeated until the leaf nodes are reached.
When \(\HSSnode\) is a leaf node, the desired matrix \(\HSSD_{\HSSnode}\) is obtained after the update, and the process terminates.
Note that the matrices \(\HSSPOST_{\HSSnode; \HSSnodech, \HSSnodech}\) in the above discussion are precisely the matrices \(\HSShD_{\HSSnodech}\) defined in~\eqref{eq:hss_construct_D_partition}.

\section{Approximate Inversion of FIOs} \label{sec:approximate_inversion_fio}

\subsection{An Inversion Algorithm Based on BF and HSS Approximations} \label{subsec:inversion_algorithm}

We now consider the problem of inverting the discrete FIO matrix \(\fiomat\) in~\eqref{eq:fio_mat}.
Suppose that \(\fiomat\) is nonsingular.
Then its inverse can be expressed as
\begin{equation} \label{eq:fio_inv}
    \fiomat^{-1} = \bigl(\fiomat^{\herm} \fiomat\bigr)^{-1} \fiomat^{\herm}.
\end{equation}

In~\cite{Feliu_Ying_2021}, the authors showed that the matrix \(\fiomat^{\herm} \fiomat\) can be approximated by an \(\hierarchical\)-matrix.
Our numerical experiments indicate that \(\fiomat^{\herm} \fiomat\) can be further compressed using an HSS matrix.
Motivated by this observation, equation~\eqref{eq:fio_inv} suggests a method for approximating \(\fiomat^{-1}\).
First, the matrix \(\fiomat\) is approximated by a BF \(\fiobf \approx \fiomat\).
Then, an HSS approximation \(\fiohss\) of \(\fioadjfio = \fiomat^{\herm} \fiomat\) is constructed by the black-box construction algorithm in Section~\ref{sec:hss_blackbox_construction}, where the BF approximation \(\fiobf^{\herm} \fiobf\) is used to perform fast matrix-vector products.
We assume that the BF approximation \(\fiobf\) is nonsingular and that the HSS approximation \(\fiohss\) is HPD.
Finally, the ULV factorization of \(\fiohss\) is computed, yielding an approximation \(\fiofactor \approx \fioadjfio^{-1}\).
Combining these approximations gives
\begin{equation} \label{eq:fio_inv_approx}
  \fiomat^{-1} \approx \fiofactor \fiobf^{\herm}.
\end{equation}
The algorithm is summarized in Algorithm~\ref{alg:fio_inv}.
\begin{algorithm}[htbp]
  \caption{Approximate inversion of the FIO} \label{alg:fio_inv}
  \begin{algorithmic}[1]
    \REQUIRE{Space and frequency grids \(\spacegrid\) and \(\freqgrid\), amplitude function \(\amp(\spacevar, \freqvar)\), phase function \(\phase(\spacevar, \freqvar)\), and accuracy parameter \(\varepsilon\).}
    \ENSURE{Approximate inverse \(\fiomat^{-1} \approx \fiofactor \fiobf^{\herm}\).}
    \STATE{Construct the BF \(\fiobf\) of \(\fiomat\).}
    \STATE{Construct the HSS approximation \(\fiohss\) using the black-box construction algorithm in Section~\ref{sec:hss_blackbox_construction}, where the BF approximation \(\fiobf^{\herm} \fiobf\) is used to perform fast matrix-vector products.}
    \STATE{Compute the ULV factorization \(\fiofactor\) of the HSS matrix \(\fiohss\).}
  \end{algorithmic}
\end{algorithm}

The approximation~\eqref{eq:fio_inv_approx} can be used as a direct solver for the linear system~\eqref{eq:fio_mat}.
The HSS inverse factor \(\fiofactor\) can also serve as a preconditioner for the approximate normal equation.
More specifically, the approximate normal equation associated with the linear system~\eqref{eq:fio_mat} is given by
\begin{equation} \label{eq:fio_approx_normal_eq}
    \fiobf^{\herm} \fiobf \fiorightmat = \fiobf^{\herm} \fioleftmat.
\end{equation}
The ULV factorization \(\fiofactor\) is then used as a preconditioner for the conjugate gradient~(CG) method applied to the HPD matrix \(\fiobf^{\herm} \fiobf\).

\subsection{Complexity Analysis} \label{subsec:complexity_analysis}

For both 1D and 2D problems, the BF approximation \(\fiobf\) can be constructed in \(\bigO(\numtot \mylog{\numtot})\) time~\cite{Li_Yang_2017}, and both \(\fiobf\) and \(\fiobf^{\herm}\) can be applied to a vector in \(\bigO(\numtot \mylog{\numtot})\) time.
It remains to analyze the complexity of the black-box construction of the HSS matrix, the ULV factorization, and the solution of the HSS system.

For 1D problems, let \(\numtot = 2^{\maxlevel} \numtot_{\maxlevel}\) be the matrix size, where \(\maxlevel\) is the maximum level of the HSS tree and \(\numtot_{\maxlevel}\) is the leaf size, which is independent of \(\numtot\).
For 2D problems, let \(\numdir = 2^{\maxlevel} \numdir_{\maxlevel}\), where \(\maxlevel\) is the maximum level of the HSS tree and \(\numdir_{\maxlevel}\) is independent of \(\numtot\).
Then \(\numtot = \numdir^{2} = 4^{\maxlevel} \numtot_{\maxlevel}\) is the matrix size, where \(\numtot_{\maxlevel} = \numdir_{\maxlevel}^{2}\) is the leaf size.
For both cases, suppose that the rank of the HSS blocks at level \(1 \leq \mylevel \leq \maxlevel\) is \(\HSSrank_{\mylevel}\), with \(\HSSrank_{0} = 0\).
The effective leaf size is defined as follows.
At the leaf level, \(\HSSleafsize_{\maxlevel} = \numtot_{\maxlevel}\).
For \(0 \leq \mylevel \leq \maxlevel - 1\), we set \(\HSSleafsize_{\mylevel} = 2 \HSSrank_{\mylevel + 1}\) for 1D problems and \(\HSSleafsize_{\mylevel} = 4 \HSSrank_{\mylevel + 1}\) for 2D problems.
The number of samples at level \(\mylevel\) is
\(\HSSsample_{\mylevel} = \HSSrank_{\mylevel} + \HSSleafsize_{\mylevel} + \oversampling\), where the oversampling parameter \(\oversampling\) is a small constant, e.g., \(\oversampling = 5\).
We assume that \(\HSSrank_{\mylevel} = \bigO(1)\) for 1D problems and \(\HSSrank_{\mylevel} = \bigO(\numdir / 2^{\mylevel})\) for 2D problems.
Under this assumption, \(\HSSleafsize_{\mylevel} = \bigO(1)\) and \(\HSSsample_{\mylevel} = \bigO(1)\) for 1D problems, while \(\HSSleafsize_{\mylevel} = \bigO(\numdir / 2^{\mylevel})\) and \(\HSSsample_{\mylevel} = \bigO(\numdir / 2^{\mylevel})\) for 2D problems.

The complexity of the black-box construction algorithm in Section~\ref{sec:hss_blackbox_construction} can be analyzed as follows.
We first consider the 1D problem.
For level \(0 \leq \mylevel \leq \maxlevel\), the cost of applying \(\HSSH^{[\mylevel]}\) to the test matrix \(\HSSOM^{[\mylevel]}\) is
\begin{equation*}
  \bigO\Biggl(
    \biggl(T_{\BF}
    + 
    \sum_{\mysublevel = \mylevel + 1}^{\maxlevel} 2^{\mysublevel} \HSSleafsize_{\mysublevel} \HSSrank_{\mysublevel}
    \biggr) \HSSsample_{\mylevel}
    \Biggr),
\end{equation*}
where \(T_{\BF}\) is the cost of applying the BF \(\fiobf\) or \(\fiobf^{\herm}\) to a vector.
For \(1 \leq \mylevel \leq \maxlevel\) and a node \(\HSSnode\) at level \(\mylevel\), the cost of computing \(\HSSU_{\HSSnode}\) is \(\bigO(\HSSsample_{\mylevel} \HSSleafsize_{\mylevel}^{2} + \HSSsample_{\mylevel} \HSSleafsize_{\mylevel} \HSSrank_{\mylevel} + \HSSleafsize_{\mylevel} \HSSrank_{\mylevel}^{2})\), and the cost of computing \(\HSScD_{\HSSnode}\) is
\(\bigO(\HSSsample_{\mylevel} \HSSleafsize_{\mylevel}^{2} + \HSSleafsize_{\mylevel}^{2} \HSSrank_{\mylevel})\).
For the root node, the cost of computing \(\HSSD_{\rootnode}\) is \(\bigO(\HSSsample_{0} \HSSleafsize_{0}^{2})\).
The postprocessing step in Section~\ref{subsec:hss_construction_postprocessing} has cost
\(\bigO\bigl(
\sum_{\mylevel = 1}^{\maxlevel} 2^{\mylevel}
(\HSSleafsize_{\mylevel} \HSSrank_{\mylevel}^{2}
+ \HSSleafsize_{\mylevel}^{2} \HSSrank_{\mylevel})
\bigr)\).
Therefore, for 1D problems, the total cost is
\begin{equation*}
  \begin{aligned}
    t_{\construct \HSS; 1 \myDim}
    & = \bigO\Biggl(
      \sum_{\mylevel = 0}^{\maxlevel} 
      \biggl(T_{\BF}
        + 
        \sum_{\mysublevel = \mylevel + 1}^{\maxlevel} 2^{\mysublevel} \HSSleafsize_{\mysublevel} \HSSrank_{\mysublevel}
        \biggr) \HSSsample_{\mylevel}
      + \sum_{\mylevel = 0}^{\maxlevel} 2^{\mylevel} 
        \bigl(\HSSsample_{\mylevel} \HSSleafsize_{\mylevel}^{2} + \HSSsample_{\mylevel} \HSSleafsize_{\mylevel} \HSSrank_{\mylevel} + \HSSleafsize_{\mylevel} \HSSrank_{\mylevel}^{2} + \HSSleafsize_{\mylevel}^{2} \HSSrank_{\mylevel}\bigr)
      \Biggr) \\
    & = \bigO\biggl(
      \sum_{\mylevel = 0}^{\maxlevel} \Bigl(\numtot \mylog{\numtot} + 2^{\maxlevel}\Bigr)
      \biggr)
    = \bigO(\numtot \mylog[2]{\numtot}).
  \end{aligned}
\end{equation*}
Similarly, for 2D problems, the total cost is
\begin{equation*}
  \begin{aligned}
    t_{\construct \HSS; 2 \myDim}
    & = \bigO\Biggl(
      \sum_{\mylevel = 0}^{\maxlevel} 
      \biggl(T_{\BF}
        + 
        \sum_{\mysublevel = \mylevel + 1}^{\maxlevel} 4^{\mysublevel} \HSSleafsize_{\mysublevel} \HSSrank_{\mysublevel}
        \biggr) \HSSsample_{\mylevel}
      + \sum_{\mylevel = 0}^{\maxlevel} 4^{\mylevel} 
        \bigl(\HSSsample_{\mylevel} \HSSleafsize_{\mylevel}^{2} + \HSSsample_{\mylevel} \HSSleafsize_{\mylevel} \HSSrank_{\mylevel} + \HSSleafsize_{\mylevel} \HSSrank_{\mylevel}^{2} + \HSSleafsize_{\mylevel}^{2} \HSSrank_{\mylevel}\bigr)
      \Biggr) \\
    & = \bigO\Biggl(
      \sum_{\mylevel = 0}^{\maxlevel}
      \biggl(
      \Bigl(\numtot \mylog{\numtot} + \numdir^{2} (\maxlevel - \mylevel)\Bigr) \frac{\numdir}{2^{\mylevel}}
       \biggr)
       + \sum_{\mylevel = 0}^{\maxlevel} 4^{\mylevel} \biggl(\frac{\numdir}{2^{\mylevel}}\biggr)^{3}
      \Biggr)
    = \bigO(\numtot^{1.5} \mylog{\numtot}).
  \end{aligned}
\end{equation*}

The ULV factorization and storage complexities follow from the rank-dependent analysis in~\cite[Theorem~6.1]{Xia_2012}.
For 1D problems, the bounded-rank assumption corresponds to case~1 with \(p=0\), which gives \(t_{\factor \HSS; 1 \myDim} = \bigO(\numtot)\) and an \(\bigO(\numtot)\) storage requirement for the ULV factors.
For 2D problems, adapting the levelwise argument of case~2(c) from the binary tree considered there to the quadtree used here gives \(t_{\factor \HSS; 2 \myDim} = \bigO(\numtot^{1.5})\) and an \(\bigO(\numtot \mylog{\numtot})\) storage requirement.
Since each stored local factor is applied only a constant number of times during an HSS solve, the solve cost has the same asymptotic order as the corresponding storage requirement.
Therefore, \(t_{\solve \HSS; 1 \myDim} = \bigO(\numtot)\) and \(t_{\solve \HSS; 2 \myDim} = \bigO(\numtot \mylog{\numtot})\).
These results are summarized in Table~\ref{tab:complexity}.
\begin{table}[htbp]
  \centering
  \begin{tabular}{c | ccccc}
  \toprule
  Dimension & \(\construct \BF\) & \(\apply \BF\) & \(\construct \HSS\) & \(\factor \HSS\) & \(\solve \HSS\) \\
  \midrule
  1D & \(\bigO(\numtot \mylog{\numtot})\) & \(\bigO(\numtot \mylog{\numtot})\) & \(\bigO(\numtot \mylog[2]{\numtot})\) & \(\bigO(\numtot)\) & \(\bigO(\numtot)\) \\
  \midrule
  2D & \(\bigO(\numtot \mylog{\numtot})\) & \(\bigO(\numtot \mylog{\numtot})\) & \(\bigO(\numtot^{1.5} \mylog{\numtot})\)  & \(\bigO(\numtot^{1.5})\) & \(\bigO(\numtot \mylog{\numtot})\) \\
  \bottomrule
  \end{tabular}
  \caption{Complexity of the main steps in the proposed approximate inversion algorithm.}
  \label{tab:complexity}
\end{table}

\section{Numerical Results} \label{sec:numerical_results}

In this section, we present numerical results for the proposed algorithm for approximating the inverse of FIOs.
The phase functions in the FIOs take the form \(\phase(\spacevar, \freqvar) = \spacevar \cdot \freqvar + \phasenonlin(\spacevar, \freqvar)\), where \(\phasenonlin(\spacevar, \freqvar)\) is periodic in \(\spacevar\).
We solve the linear system~\eqref{eq:fio_mat} through the approximate normal equation~\eqref{eq:fio_approx_normal_eq}.
For the direct solver, the solution is computed by solving the HSS system using the ULV factorization \(\fiofactor\).
For the iterative solver, we use the MATLAB command \texttt{pcg}, and the solution is computed by the conjugate gradient~(CG) method either without preconditioning or with the ULV factorization \(\fiofactor\) as a preconditioner.
For 1D problems, we let \(\numtot = \numdir\) range over \(\{2^{10}, 2^{12}, 2^{14}, 2^{16}, 2^{18}\}\).
For 2D problems, we let \(\numdir\) range over \(\{64, 128, 256, 512\}\) and set \(\numtot = \numdir^{2}\).
The real and imaginary parts of \(\fiorightmat\) have independent standard Gaussian entries, and the right-hand side is set to
\(\fioleftmat = \fiobf \fiorightmat\).
This choice assesses the HSS-based inversion separately from the BF approximation error, which is reported independently.
All algorithms are implemented in MATLAB R2023b, and the BF is computed using the code in~\cite{Li_2016_code}.
All experiments are carried out on a server with an Intel Gold 6226R CPU at 2.90 GHz and 1000.6 GB of RAM.

The following quantities are used to evaluate the performance of the algorithms.
We denote the walltimes for BF construction, HSS construction, ULV factorization, and HSS solution by \(t_{\construct \BF}\), \(t_{\construct \HSS}\), \(t_{\factor \HSS}\), and \(t_{\solve \HSS}\), respectively.
All walltimes are reported in seconds.
The total offline setup cost is \(t_{\pre} = t_{\construct \BF} + t_{\construct \HSS} + t_{\factor \HSS}\), while \(t_{\solve \HSS}\) measures the online HSS solve time.
The tolerance parameters for the BF construction and HSS construction are fixed at \(10^{-8}\) and \(10^{-3}\), respectively.
The relative approximation errors are defined as
\begin{equation*}
  e_{\construct \BF}
  =
  \frac{\|\fiobf \fiorightmat - \fiomat \fiorightmat\|_{2}}
       {\|\fiomat \fiorightmat\|_{2}},
  \quad
  e_{\construct \HSS}
  =
  \frac{\|\fiohss \fiorightmat - \fiobf^{\herm} \fiobf \fiorightmat\|_{2}}
       {\|\fiohss \fiorightmat\|_{2}}.
\end{equation*}
We estimate \(e_{\construct \BF}\) by restricting both norms to \(256\) randomly sampled output entries, without explicitly forming \(\fiomat\).
For a computed solution \(\fiorightmatsol\), the relative error is
\(e_{\solve} = \|\fiorightmat - \fiorightmatsol\|_{2} / \|\fiorightmat\|_{2}\),
with \(e_{\direct}\) denoting the direct solver error.
For CG with or without preconditioning, the relative residual tolerance is \(10^{-12}\).
For CG applied to~\eqref{eq:fio_approx_normal_eq}, the walltime, iteration count, and solution error are denoted by
\((t_{\iter}, n_{\iter}, e_{\iter})\) without preconditioning and
\((t_{\piter}, n_{\piter}, e_{\piter})\) with preconditioning by \(\fiofactor\).

\subsection{1D FIO With Variable Amplitude} \label{subsec:1d_var_amp}

We begin with a 1D FIO whose amplitude and phase functions are given by
\begin{equation*}
  \begin{gathered}
    \amp(\spacevar, \freqvar) = \sum_{k = 1}^{m} \exp\biggl(-
    \frac{(\spacevar - \spacevar_{k})^{2} + ((\freqvar - \freqvar_{k}) / \numtot)^{2}}{\sigma^{2}}\biggr), \\
    \phase(\spacevar, \freqvar) = \spacevar \cdot \freqvar + \phasenonlin(\spacevar, \freqvar),
  \quad
  \phasenonlin(\spacevar, \freqvar) = \frac{2 + \sin(2 \pi \spacevar)}{8} |\freqvar| ,
  \end{gathered}
\end{equation*}
where \(\spacevar_{k}\) and \(\freqvar_{k}\) are randomly generated points in \([0, 1] \times [-\numtot / 2, \numtot / 2]\) for \(1 \leq k \leq m\), and \(\sigma^{2}\) is a positive constant.
In our experiments, \(m\) is set to \(10\) and \(\sigma^{2}\) is set to \(0.1\).

\begin{figure}[tbhp]
    \centering
    \setlength{\directplotheight}{0.28\textwidth}

    \includegraphics[height=\directplotheight]{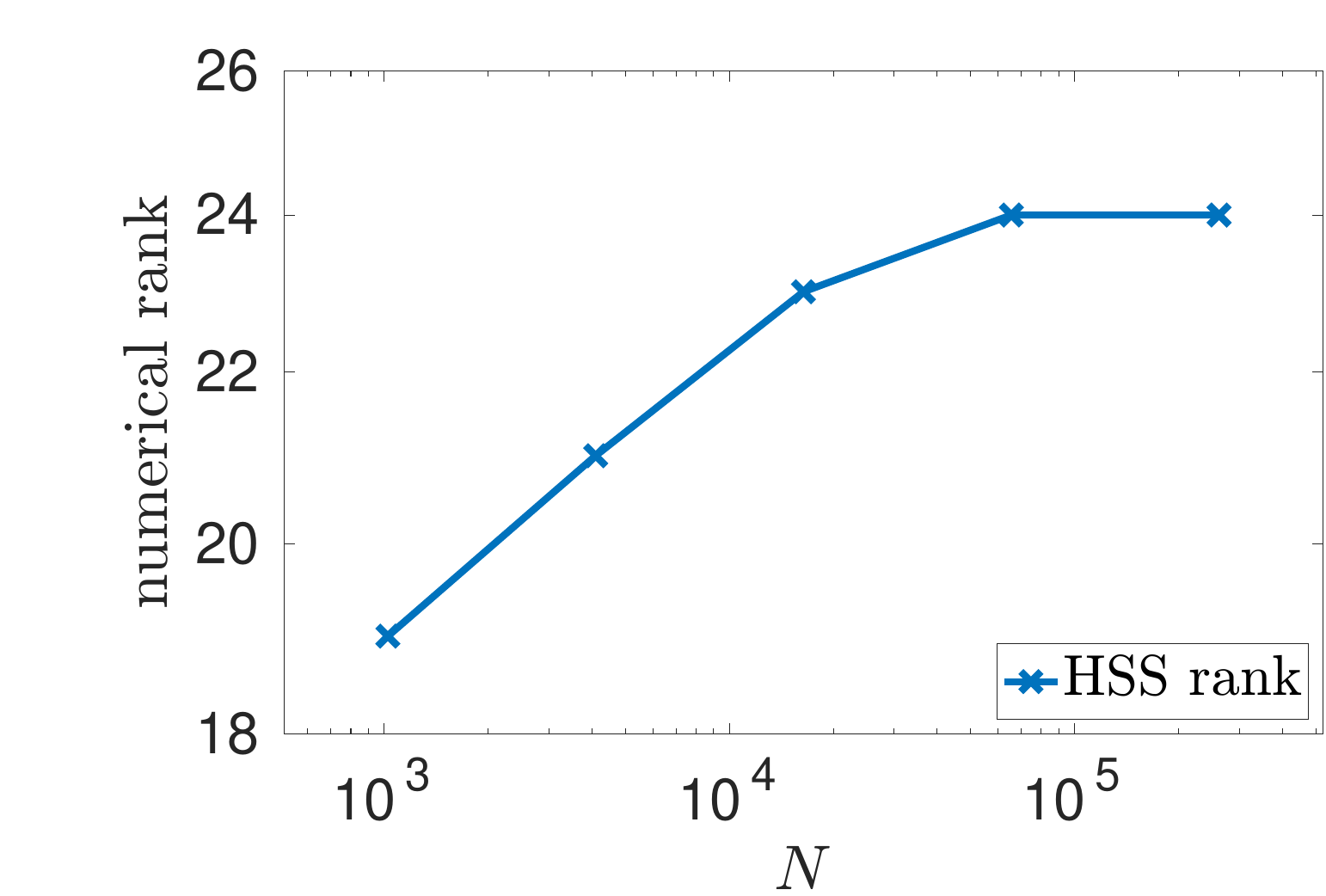}%
    \hfill
    \includegraphics[height=\directplotheight]{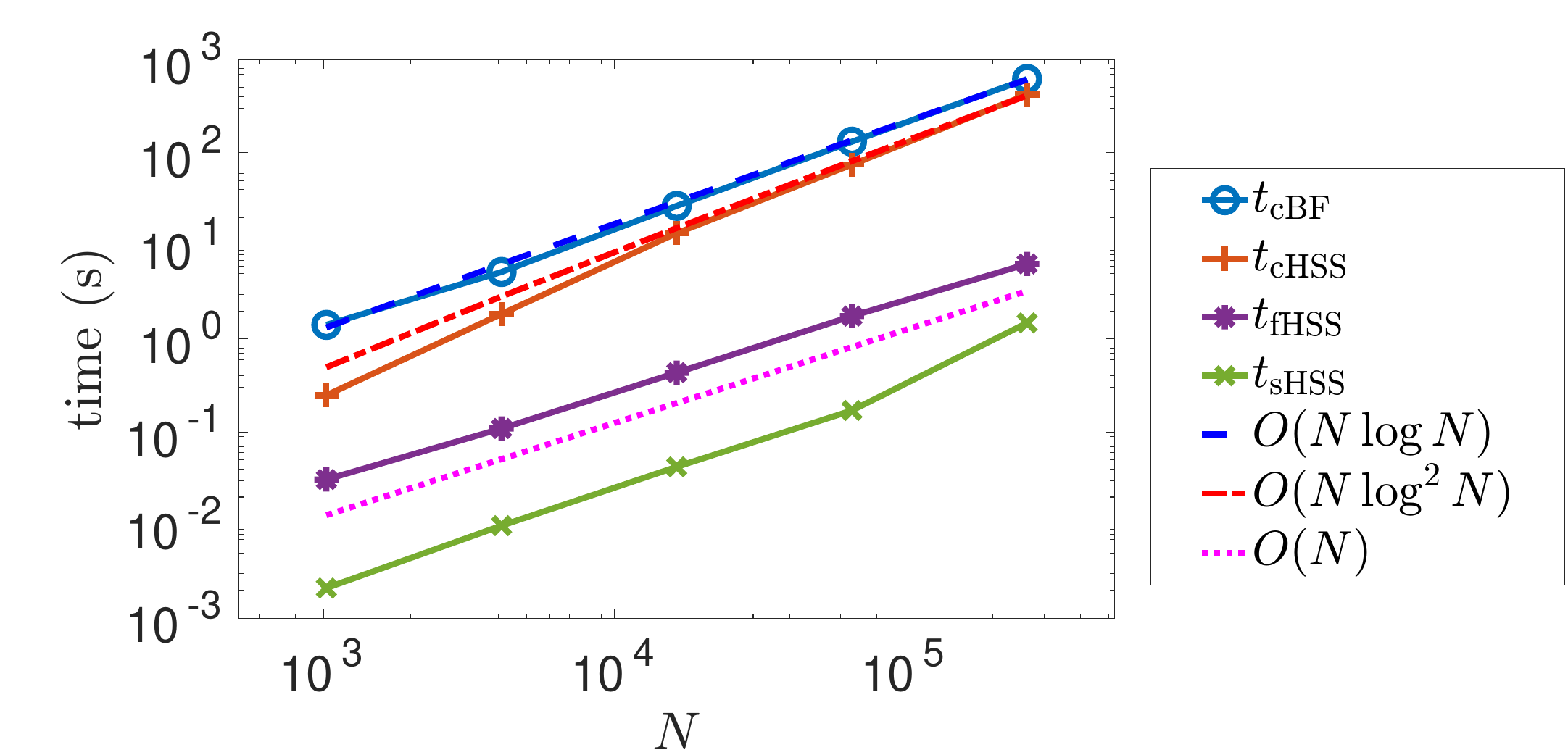}

    \caption{HSS rank and time scaling of the direct solver for a
    1D FIO with variable amplitude.}
    \label{fig:1d_var_amplitude_direct}
\end{figure}

\begin{table}[tbhp]
\centering
\begin{tabular}{c | cc | cc | cc | c}
\toprule
\(\numtot\) & \(t_{\construct \BF}\) & \(e_{\construct \BF}\) & \(t_{\construct \HSS}\) & \(e_{\construct \HSS}\) & \(t_{\factor \HSS}\) & \(t_{\solve \HSS}\) & \(e_{\direct}\) \\ 
\midrule
\(2^{10}\) & 1.4e+00 & 3.7e-06 & 2.5e-01 & 2.0e-05 & 3.1e-02 & 2.1e-03 & 1.2e-04 \\ 
\midrule
\(2^{12}\) & 5.2e+00 & 3.2e-06 & 1.9e+00 & 2.9e-05 & 1.1e-01 & 9.9e-03 & 1.7e-04 \\ 
\midrule
\(2^{14}\) & 2.7e+01 & 3.6e-06 & 1.4e+01 & 3.0e-05 & 4.3e-01 & 4.2e-02 & 1.9e-04 \\ 
\midrule
\(2^{16}\) & 1.3e+02 & 5.3e-06 & 7.4e+01 & 3.1e-05 & 1.8e+00 & 1.7e-01 & 1.9e-04 \\ 
\midrule
\(2^{18}\) & 6.2e+02 & 4.2e-06 & 4.2e+02 & 3.3e-05 & 6.4e+00 & 1.5e+00 & 2.1e-04 \\ 
\bottomrule
\end{tabular}
\caption{Results of the direct solver for a 1D FIO with variable amplitude.}
\label{tab:1d_var_amplitude_direct}
\end{table}

\begin{table}[tbhp]
\centering
\begin{tabular}{c | ccc | cccc}
\toprule
\(\numtot\) & \(t_{\iter}\) & \(n_{\iter}\) & \(e_{\iter}\) & \(t_{\pre}\) & \(t_{\piter}\) & \(n_{\piter}\) & \(e_{\piter}\) \\ 
\midrule
\(2^{10}\) & 5.0e-01 & 141 & 5.7e-12 & 1.7e+00 & 3.5e-02 & 4 & 1.5e-14 \\ 
\midrule
\(2^{12}\) & 2.9e+00 & 159 & 5.3e-12 & 7.2e+00 & 1.4e-01 & 4 & 5.1e-13 \\ 
\midrule
\(2^{14}\) & 1.5e+01 & 164 & 6.0e-12 & 4.1e+01 & 6.3e-01 & 4 & 6.8e-13 \\ 
\midrule
\(2^{16}\) & 7.3e+01 & 166 & 5.9e-12 & 2.1e+02 & 2.8e+00 & 4 & 2.2e-12 \\ 
\midrule
\(2^{18}\) & 3.2e+02 & 166 & 6.1e-12 & 1.0e+03 & 2.0e+01 & 5 & 1.0e-13 \\ 
\bottomrule
\end{tabular}
\caption{Results of the iterative solver for a 1D FIO with variable amplitude.}
\label{tab:1d_var_amplitude_iterative}
\end{table}

Figure~\ref{fig:1d_var_amplitude_direct} and Table~\ref{tab:1d_var_amplitude_direct} report the HSS ranks, walltimes, and accuracy of the direct solver.
The HSS ranks remain essentially bounded as \(\numtot\) increases, which is consistent with the rank assumption used in the complexity analysis in Section~\ref{subsec:complexity_analysis}.
The observed walltimes also agree with the predicted scaling: the BF construction and application scale as \(\bigO(\numtot \mylog{\numtot})\), while the HSS construction and the ULV factorization exhibit the expected \(\bigO(\numtot \mylog[2]{\numtot})\) and \(\bigO(\numtot)\) behavior, respectively.
In this example, the BF construction takes the largest part of the time.
The HSS approximation and the solution computed by the direct solver are stable for all tested problem sizes, with \(e_{\construct \mathrm{HSS}}\) on the order of \(10^{-5}\) and \(e_{\mathrm{\solve}}\) on the order of \(10^{-4}\).

Table~\ref{tab:1d_var_amplitude_iterative} summarizes the results of the iterative solver.
Without preconditioning, CG requires about \(160\) iterations for all values of \(\numtot\).
In contrast, preconditioning CG with \(\fiofactor\) reduces the iteration count to only \(4\) or \(5\), with solution errors below
\(3 \times 10^{-12}\) for all tested problem sizes.
Moreover, the iteration time \(t_{\iter}\) for unpreconditioned CG is nearly \(20\) times that for preconditioned CG across all tested problem sizes.
These results show that the proposed approximate inverse is an effective preconditioner for the normal equation.

\subsection{2D FIO With Constant Amplitude} \label{subsec:2d_const_amp}

In this example, we consider a constant amplitude function \(\amp(\spacevar, \freqvar) = 1\) and the phase function
\begin{equation*}
  \begin{gathered}
    \phase(\spacevar, \freqvar) = \spacevar \cdot \freqvar + \phasenonlin(\spacevar, \freqvar),
    \quad
    \phasenonlin(\spacevar, \freqvar) = \sqrt{\phasecoeff_{1}^{2}(\spacevar) \freqvarsimple_{1}^{2} + \phasecoeff_{2}^{2}(\spacevar) \freqvarsimple_{2}^{2}}, \\
    \phasecoeff_{1}(\spacevar) = \frac{2 + \sin(2 \pi \spacevarsimple_{1}) \sin(2 \pi \spacevarsimple_{2})}{16},
    \quad
    \phasecoeff_{2}(\spacevar) = \frac{2 + \cos(2 \pi \spacevarsimple_{1}) \cos(2 \pi \spacevarsimple_{2})}{16}.
  \end{gathered}
\end{equation*}

\begin{figure}[tbhp]
    \centering
    \setlength{\directplotheight}{0.27\textwidth}

    \includegraphics[height=\directplotheight]{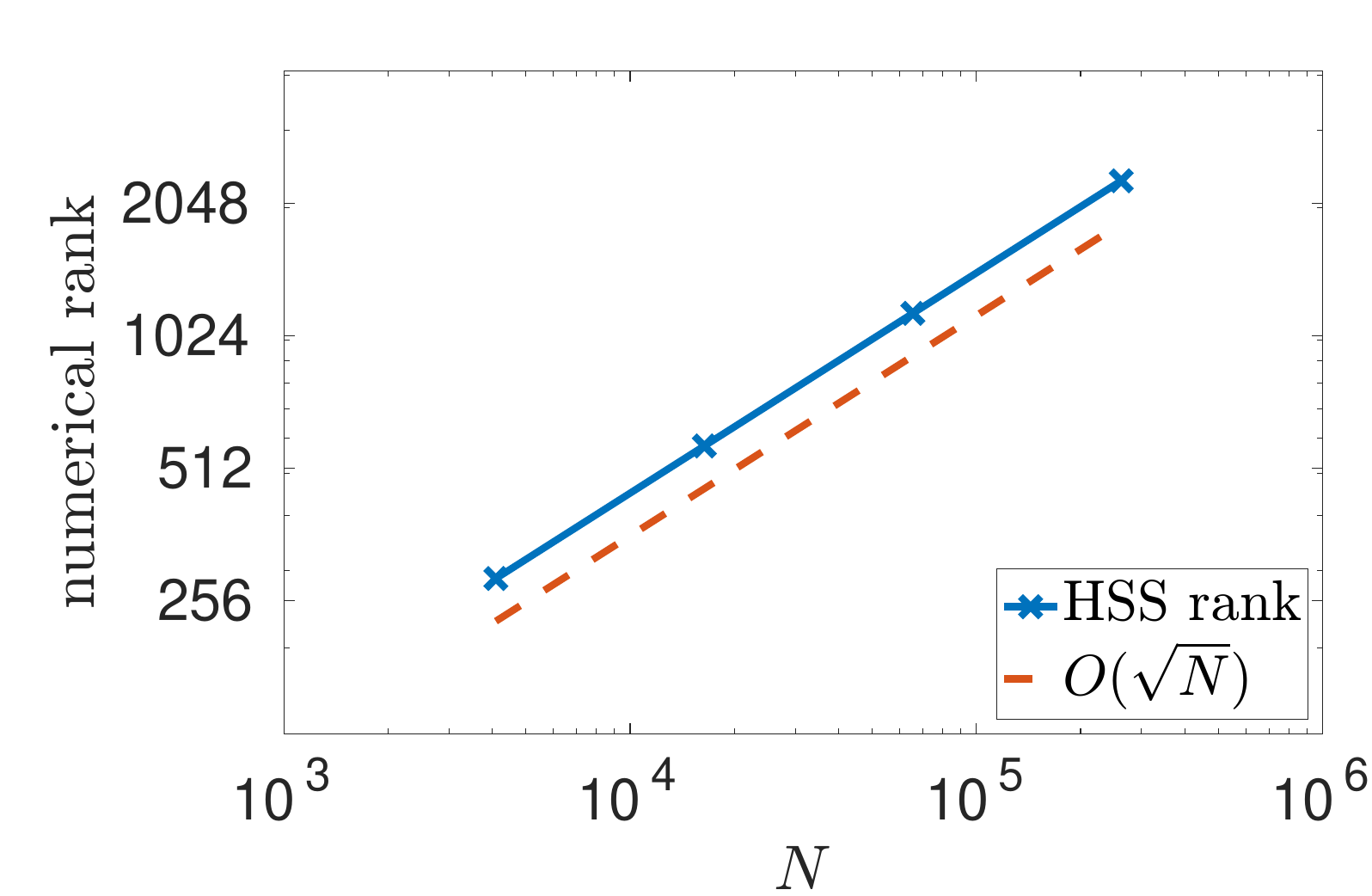}%
    \hfill
    \includegraphics[height=\directplotheight]{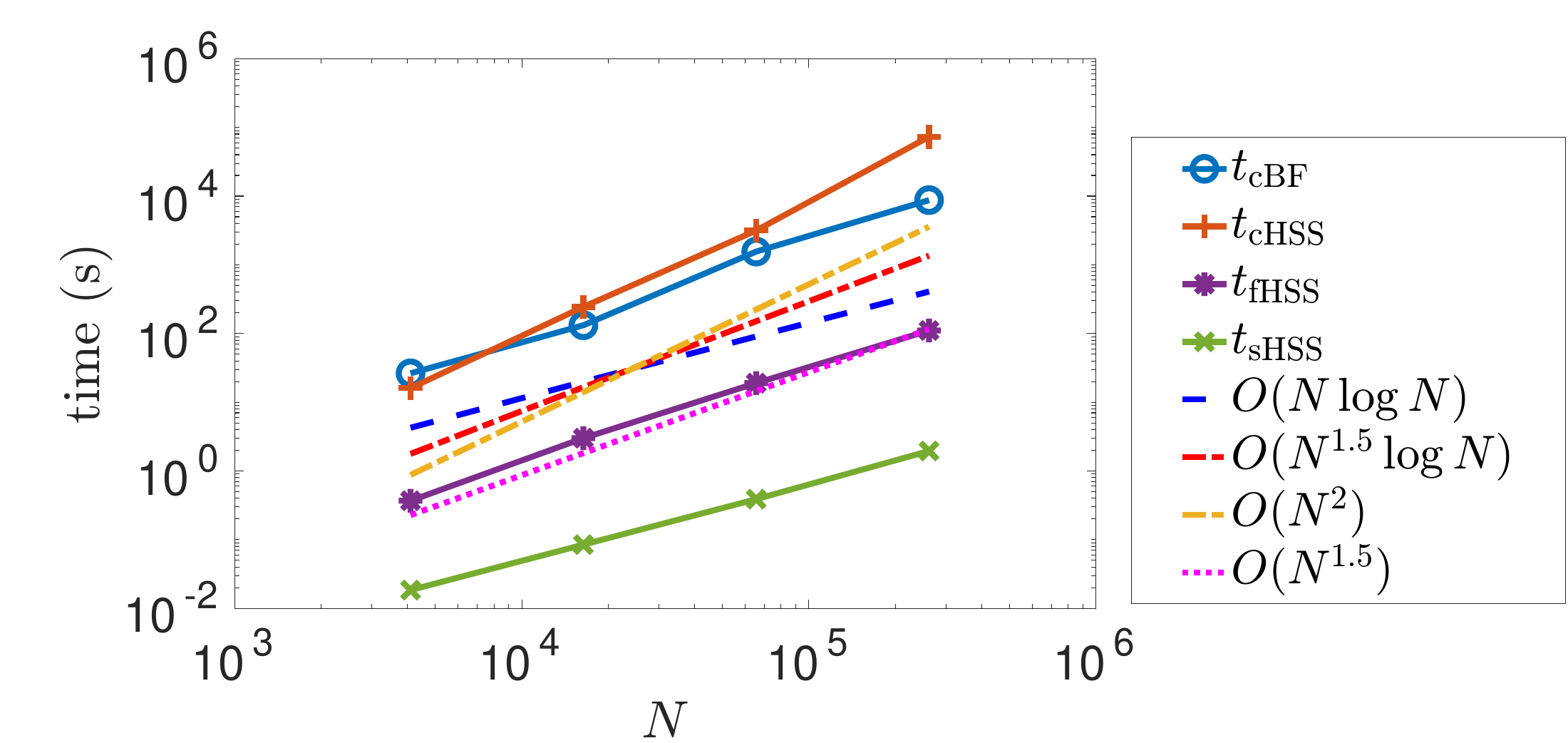}

    \caption{HSS rank and time scaling of the direct solver for a
    2D FIO with constant amplitude.}
    \label{fig:2d_const_amplitude_direct}
\end{figure}

\begin{table}[tbhp]
\centering
\begin{tabular}{c | cc | cc | cc | c}
\toprule
\(\numtot\) & \(t_{\construct \BF}\) & \(e_{\construct \BF}\) & \(t_{\construct \HSS}\) & \(e_{\construct \HSS}\) & \(t_{\factor \HSS}\) & \(t_{\solve \HSS}\) & \(e_{\direct}\) \\ 
\midrule
\(64^{2}\) & 2.6e+01 & 3.4e-05 & 1.6e+01 & 8.5e-04 & 3.7e-01 & 1.8e-02 & 9.7e-04 \\ 
\midrule
\(128^{2}\) & 1.3e+02 & 3.0e-05 & 2.4e+02 & 1.2e-03 & 3.0e+00 & 8.4e-02 & 1.4e-03 \\ 
\midrule
\(256^{2}\) & 1.6e+03 & 1.6e-05 & 3.2e+03 & 1.6e-03 & 1.9e+01 & 3.9e-01 & 1.8e-03 \\ 
\midrule
\(512^{2}\) & 8.7e+03 & 2.2e-05 & 7.3e+04 & 1.9e-03 & 1.1e+02 & 1.9e+00 & 2.1e-03 \\ 
\bottomrule
\end{tabular}
\caption{Results of the direct solver for a 2D FIO with constant amplitude.}
\label{tab:2d_const_amplitude_direct}
\end{table}

\begin{table}[tbhp]
\centering
\begin{tabular}{c | ccc | cccc}
\toprule
\(\numtot\) & \(t_{\iter}\) & \(n_{\iter}\) & \(e_{\iter}\) & \(t_{\pre}\) & \(t_{\piter}\) & \(n_{\piter}\) & \(e_{\piter}\) \\ 
\midrule
\(64^{2}\) & 3.5e+00 & 30 & 4.3e-13 & 4.3e+01 & 8.0e-01 & 5 & 7.2e-13 \\ 
\midrule
\(128^{2}\) & 3.1e+01 & 31 & 8.8e-13 & 3.8e+02 & 7.1e+00 & 6 & 2.7e-14 \\ 
\midrule
\(256^{2}\) & 2.0e+02 & 32 & 7.4e-13 & 4.7e+03 & 4.6e+01 & 6 & 3.0e-13 \\ 
\midrule
\(512^{2}\) & 1.1e+03 & 32 & 9.3e-13 & 8.2e+04 & 2.4e+02 & 6 & 7.0e-13 \\ 
\bottomrule
\end{tabular}
\caption{Results of the iterative solver for a 2D FIO with constant amplitude.}
\label{tab:2d_const_amplitude_iterative}
\end{table}

Figure~\ref{fig:2d_const_amplitude_direct} and Table~\ref{tab:2d_const_amplitude_direct} report the HSS ranks, walltimes, and accuracy of the direct solver.
The observed HSS ranks of the constructed HSS representations grow consistently with the expected \(\bigO(\sqrt{\numtot})\) scaling.
The walltimes for BF construction, HSS ULV factorization, and solving the HSS system also agree with the complexity estimates in Section~\ref{subsec:complexity_analysis}.
The HSS construction time, however, exhibits an empirical scaling closer to \(\bigO(\numtot^{2})\), rather than the predicted \(\bigO(\numtot^{1.5} \mylog{\numtot})\) scaling.
This discrepancy is caused by the practical implementation of the BF in~\eqref{eq:bf}.
In the implementation, the super-index of the outermost factors \(\BFU\), \(\BFG\), \(\BFH\), and \(\BFV\) is not \(\maxlevel\), but rather \(\max\{\maxlevel - \maxleveloffset, \midlevel\}\), where \(\maxlevel = \mylog_{2}{\numdir}\).
Consequently, the number of \(\BFG\) and \(\BFH\) factors in the BF is
\begin{equation*}
  \max\{\maxlevel - \maxleveloffset, \midlevel\} - \midlevel
  =
  \max\{\floor(\maxlevel / 2) - \maxleveloffset, 0\}.
\end{equation*}
For \(\numtot \leq 128^{2}\), we have \(\maxlevel \leq 7\), and hence no \(\BFG\) or \(\BFH\) factors are present.
For \(\numtot = 256^{2}\) and \(512^{2}\), we have \(\maxlevel = 8\) and \(9\), respectively, so only one \(\BFG\) factor and one \(\BFH\) factor are present.
In these regimes, the application cost of the BF deteriorates to \(\bigO(\numtot^{1.5})\), which in turn makes the HSS construction time scale empirically like \(\bigO(\numtot^{2})\).
Note that the HSS construction becomes the most time-consuming part in this example.
Despite this implementation-dependent slowdown, the relative errors of the HSS approximation and the direct solver remain on the order of \(10^{-3}\) for all tested problem sizes, indicating that the HSS construction and the direct solver are stable.

Table~\ref{tab:2d_const_amplitude_iterative} summarizes the results of the iterative solver.
Without preconditioning, CG requires about \(30\) iterations for all values of \(\numtot\), suggesting that the condition number of the approximate normal equation~\eqref{eq:fio_approx_normal_eq} remains moderate and does not grow significantly with \(\numtot\).
Preconditioning CG with \(\fiofactor\) further reduces the iteration
count to \(5\) or \(6\), while maintaining solution errors below
\(10^{-12}\) for all tested problem sizes.
The iteration time \(t_{\iter}\) for unpreconditioned CG is nearly \(5\) times that for preconditioned CG.
These results show that the proposed approximate inverse remains an effective preconditioner for 2D FIOs with constant amplitude.

\subsection{2D FIO With Variable Amplitude} \label{subsec:2d_var_amp}

In this example, the amplitude and phase functions are given by
\begin{equation*}
  \begin{gathered}
    \amp(\spacevar, \freqvar) = \besselh_{0}^{(1)} (2 \pi \phasenonlin(\spacevar, \freqvar)) \e^{-2 \pi \imath \phasenonlin(\spacevar, \freqvar)}, \\
    \phase(\spacevar, \freqvar) = \spacevar \cdot \freqvar + \phasenonlin(\spacevar, \freqvar),
    \quad
    \phasenonlin(\spacevar, \freqvar) = \phasecoeff(\spacevar) \|\freqvar\|,
    \quad
    \phasecoeff(\spacevar) = \frac{3 + \sin(2 \pi \spacevarsimple_{1}) \sin(2 \pi \spacevarsimple_{2})}{8},
  \end{gathered}
\end{equation*}
where \(\besselh_{0}^{(1)}\) denotes the Hankel function of the first kind of order \(0\), and \(\|\cdot\|\) denotes the Euclidean norm.

\begin{figure}[tbhp]
    \centering
    \setlength{\directplotheight}{0.27\textwidth}

    \includegraphics[height=\directplotheight]{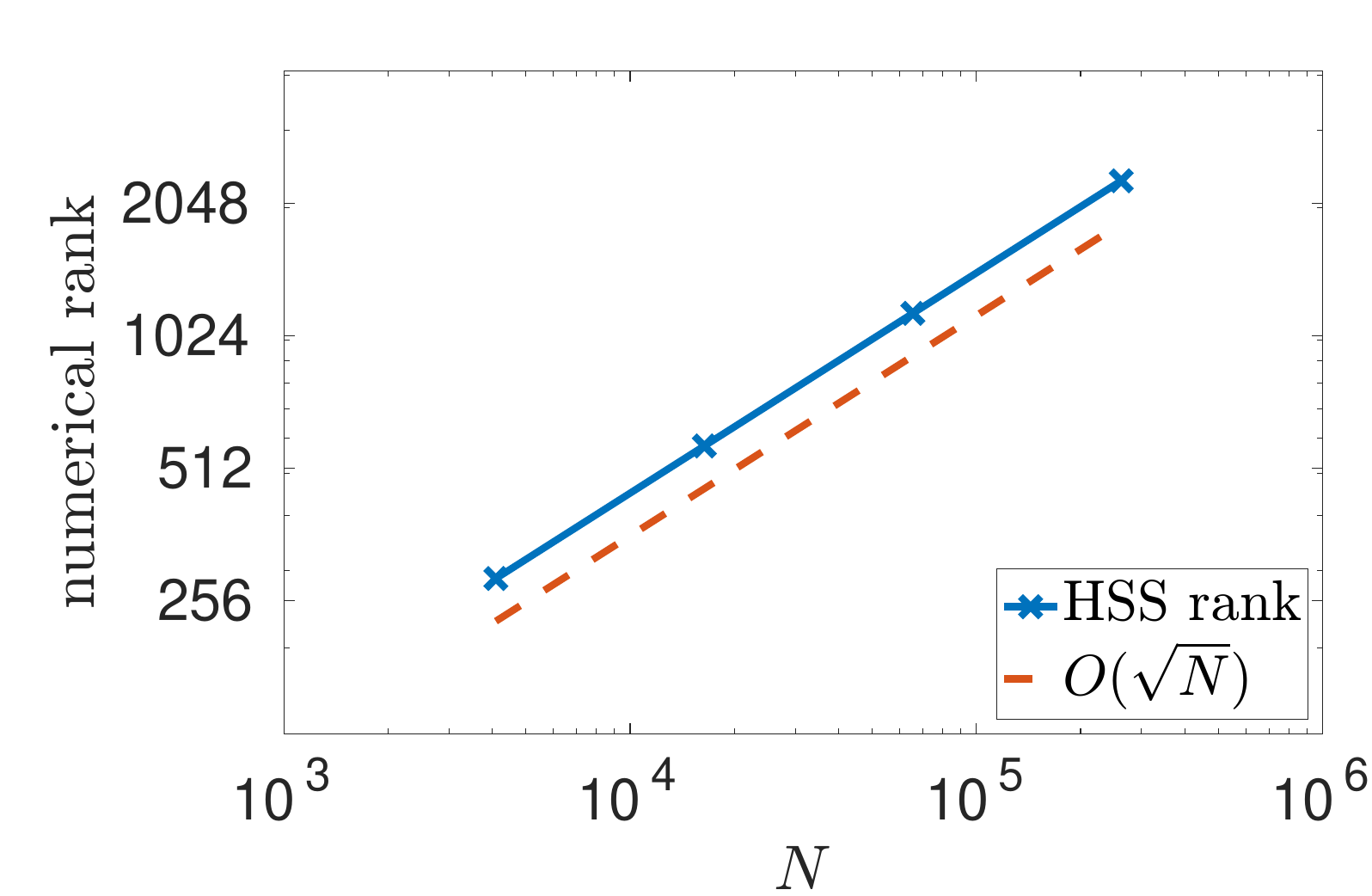}%
    \hfill
    \includegraphics[height=\directplotheight]{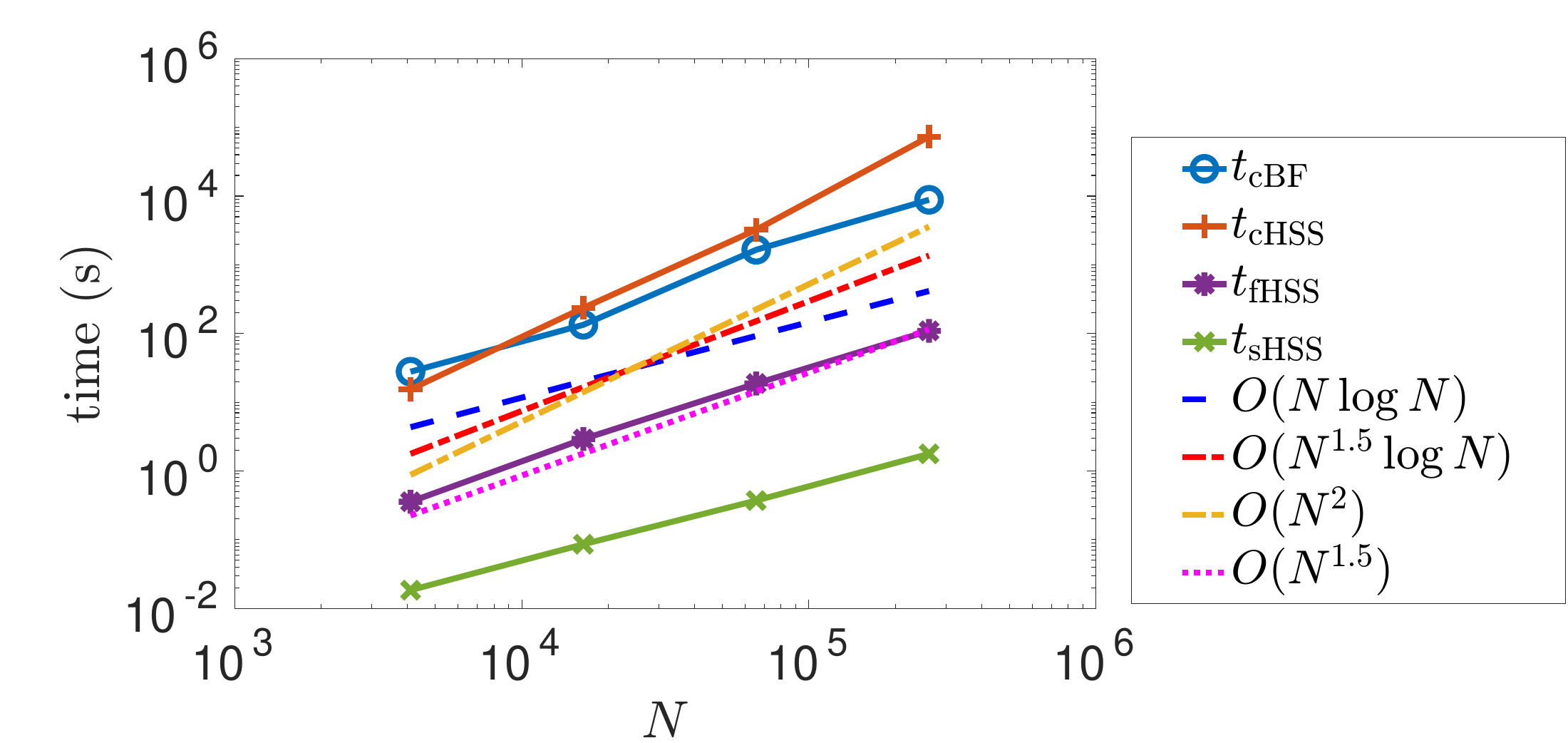}

    \caption{HSS rank and time scaling of the direct solver for a
    2D FIO with variable amplitude.}
    \label{fig:2d_var_amplitude_direct}
\end{figure}

\begin{table}[tbhp]
\centering
\begin{tabular}{c | cc | cc | cc | c}
\toprule
\(\numtot\) & \(t_{\construct \BF}\) & \(e_{\construct \BF}\) & \(t_{\construct \HSS}\) & \(e_{\construct \HSS}\) & \(t_{\factor \HSS}\) & \(t_{\solve \HSS}\) & \(e_{\direct}\) \\ 
\midrule
\(64^{2}\) & 2.8e+01 & 2.1e-03 & 1.5e+01 & 1.3e-03 & 3.6e-01 & 1.8e-02 & 3.8e-03 \\ 
\midrule
\(128^{2}\) & 1.3e+02 & 5.3e-04 & 2.4e+02 & 1.6e-03 & 2.9e+00 & 8.6e-02 & 3.9e-03 \\ 
\midrule
\(256^{2}\) & 1.7e+03 & 2.6e-04 & 3.2e+03 & 1.6e-03 & 1.9e+01 & 3.7e-01 & 3.4e-03 \\ 
\midrule
\(512^{2}\) & 8.8e+03 & 2.3e-04 & 7.3e+04 & 1.8e-03 & 1.1e+02 & 1.8e+00 & 3.2e-03 \\ 
\bottomrule
\end{tabular}
\caption{Results of the direct solver for a 2D FIO with variable amplitude.}
\label{tab:2d_var_amplitude_direct}
\end{table}

\begin{table}[tbhp]
\centering
\begin{tabular}{c | ccc | cccc}
\toprule
\(\numtot\) & \(t_{\iter}\) & \(n_{\iter}\) & \(e_{\iter}\) & \(t_{\pre}\) & \(t_{\piter}\) & \(n_{\piter}\) & \(e_{\piter}\) \\ 
\midrule
\(64^{2}\) & 1.3e+01 & 114 & 2.0e-12 & 4.4e+01 & 8.9e-01 & 6 & 1.1e-12 \\ 
\midrule
\(128^{2}\) & 1.6e+02 & 161 & 1.5e-12 & 3.7e+02 & 7.2e+00 & 6 & 1.4e-12 \\ 
\midrule
\(256^{2}\) & 1.3e+03 & 223 & 1.8e-12 & 4.9e+03 & 4.8e+01 & 6 & 8.7e-13 \\ 
\midrule
\(512^{2}\) & 1.1e+04 & 317 & 1.4e-12 & 8.2e+04 & 2.5e+02 & 6 & 6.8e-13 \\ 
\bottomrule
\end{tabular}
\caption{Results of the iterative solver for a 2D FIO with variable amplitude.}
\label{tab:2d_var_amplitude_iterative}
\end{table}

Figure~\ref{fig:2d_var_amplitude_direct} and Table~\ref{tab:2d_var_amplitude_direct} report the HSS ranks, walltimes, and accuracy of the direct solver.
The observed HSS ranks of the constructed HSS representations grow consistently with the expected \(\bigO(\sqrt{\numtot})\) scaling, as in the constant-amplitude case.
The walltimes of all components of the direct solver, except for HSS construction, agree with the complexity estimates in Section~\ref{subsec:complexity_analysis}.
As in Section~\ref{subsec:2d_const_amp}, the empirical \(\bigO(\numtot^{2})\) scaling of the HSS construction time is caused by the implementation-dependent cost of applying the BF.
The most time-consuming part of the direct solver is again the HSS construction.
The relative errors of the HSS approximation and the solution computed by the direct solver remain on the order of \(10^{-3}\) for all tested problem sizes, indicating that the HSS construction and the direct solver are stable.

Table~\ref{tab:2d_var_amplitude_iterative} summarizes the results of the iterative solver.
Without preconditioning, the number of CG iterations increases from \(114\) to \(317\) as \(\numtot\) increases, suggesting that the conditioning of the approximate normal equation~\eqref{eq:fio_approx_normal_eq} deteriorates with the problem size.
In contrast, using \(\fiofactor\) as a preconditioner reduces the iteration count to \(6\) for all tested problem sizes, while maintaining a solution error on the order of \(10^{-12}\).
For \(\numtot = 256^{2}\) and \(512^{2}\), the iteration time \(t_{\iter}\) for unpreconditioned CG is nearly \(27\) and \(44\) times that for preconditioned CG, respectively.
These results demonstrate that the proposed preconditioner is effective for this problem, and the speedup becomes more significant as \(\numtot\) increases.

\section{Conclusions} \label{sec:conclusions}

In this paper, we introduced an algorithm for approximating the inverse of the discrete FIO \(\fiomat\).
The offline stage of the algorithm consists of three steps: constructing a BF approximation \(\fiobf \approx \fiomat\), constructing an HSS approximation \(\fiohss \approx \fioadjfio = \fiomat^{\herm} \fiomat\), and computing the ULV factorization of \(\fiohss\) to obtain \(\fiofactor \approx \fioadjfio^{-1}\).
Using these offline components, the approximate inverse of the discrete FIO is applied in the online stage as \(\fiomat^{-1} \approx \fiofactor \fiobf^{\herm}\).
For the HSS construction, we developed a black-box method that only requires matrix-vector products with the target matrix.
This method is efficient in the present setting because these matrix-vector products can be performed rapidly using the BF approximation.

The proposed algorithm provides an HSS-based alternative to the \(\hierarchical\)-matrix and HIF approach in~\cite{Feliu_Ying_2021}.
It combines the nested shared bases of the HSS representation with a black-box construction algorithm and ULV factorization.
The HSS representation reuses basis information across levels, and its generators support a recursive ULV factorization.
The black-box construction uses independent random test matrices at different levels, allowing the sampling to be adapted to the numerical ranks at each level.

For 1D FIOs, the offline costs for constructing the BF approximation, constructing the HSS approximation, and computing the ULV factorization are \(\bigO(\numtot \mylog{\numtot})\), \(\bigO(\numtot \mylog[2]{\numtot})\), and \(\bigO(\numtot)\), respectively, while the online solution cost is \(\bigO(\numtot \mylog{\numtot})\).
For 2D FIOs, the corresponding offline costs are \(\bigO(\numtot \mylog{\numtot})\), \(\bigO(\numtot^{1.5} \mylog{\numtot})\), and \(\bigO(\numtot^{1.5})\), respectively, while the online solution cost is \(\bigO(\numtot \mylog{\numtot})\).
The resulting approximate inverse can be used as a direct solver, while \(\fiofactor\) can be used as a preconditioner for CG applied to the approximate normal equation.
In the reported experiments, this preconditioner kept the iteration count small and nearly independent of the problem size, while
substantially reducing the iteration time relative to unpreconditioned CG.

There are several directions for future work.
First, the proposed algorithm can be extended to higher-dimensional FIOs in a straightforward manner by using the same framework.
However, in higher dimensions, the HSS ranks may grow with the problem size, leading to higher computational complexity.
Second, a parallel implementation of the proposed algorithm could be developed to further reduce the walltime, especially for large-scale problems.
Indeed, the main components of the algorithm, including BF construction, HSS construction, and ULV factorization, are all well suited for parallelization.
Third, establishing the HSS property of the matrix \(\fioadjfio\) and analyzing the corresponding HSS ranks are important theoretical questions.
Such results would provide a rigorous foundation for the complexity analysis of the proposed algorithm.
Finally, in~\cite{Feliu_Ying_2021}, the authors used \(\hierarchical\)-matrices with strong admissibility to approximate
\(\fioadjfio=\fiomat^{\herm}\fiomat\) for 2D problems.
Their experiments indicated that the ranks of the admissible blocks grow approximately as \(\bigO(\mylog{\numtot})\).
It would therefore be interesting to investigate whether \(\hierarchical^{2}\)-matrix approximations~\cite{Hackbusch_Khoromskij_Sauter_2000} and their inverses~\cite{Boukaram_Keyes_Li_Liu_Turkiyyah_2026, Ma_Deshmukh_Yokota_2022} can be used for high-dimensional problems.
This may lead to a more efficient algorithm for approximating the inverse of FIOs in higher dimensions.

\bibliographystyle{siam}
\bibliography{ref} 
 
\end{document}